\documentclass[12pt]{article}
\usepackage{amsfonts,amssymb,amscd}
\usepackage{oldgerm}
\usepackage{fullpage}
\newtheorem{teo}{Theorem}[section]

\newtheorem{lema}[teo]{Lemma}

\newtheorem{defnc}[teo]{Definition}

\newtheorem{prop}[teo]{Proposition}

\newtheorem{obs}[teo]{Remark}

\newtheorem{coro}[teo]{Corollary}

\newcommand{\C}{{\mathbb C}}

\newcommand{\R}{{\mathbb R}}

\newcommand{\Z}{{\mathbb Z}}

\newcommand{\N}{{\mathbb N}}

\newcommand{\fii}{{\Phi_{sg}}}

\newcommand{\parlx}{{\partial /\partial x}}

\newcommand{\parly}{{\partial /\partial y}}

\newcommand{\partx}{{\frac{\partial}{\partial x}}}

\newcommand{\party}{{\frac{\partial}{\partial y}}}

\newcommand{\fol}{{\mathcal F}}

\newcommand{\cale}{{\mathcal E}}

\newcommand{\cals}{{\mathcal S}}

\newcommand{\Pp}{{\mathbb P}}

\newcommand{\pp}{{\mathbb C} {\mathbb P} (2)}

\newcommand{\wifol}{\widetilde{\mathcal F}}

\newcommand{\widxx}{\widetilde{X}}

\newcommand{\wiX}{\widetilde{Y}}

\newcommand{\wicc}{\widetilde{\mathbb C}^2}

\newcommand{\qed}{\hfill \fbox{}}

\begin{document}

\title{Complete polynomial vector fields on $\C^2$, {\sc Part I}}

\author{Julio C. Rebelo}

\date{}

\maketitle
\thispagestyle{empty} \def\IMSmarkvadjust{0 pt}
\def\IMSmarkhadjust{0 pt}
\def\IMSmarkhpadding{0 pt}
\def\IMSpubltext{Published in modified form:}
\def\SBIMSMark#1#2#3{
 \font\SBF=cmss10 at 10 true pt
 \font\SBI=cmssi10 at 10 true pt
 \setbox0=\hbox{\SBF \hbox to \IMSmarkhpadding{\relax}
                Stony Brook IMS Preprint \##1}
 \setbox2=\hbox to \wd0{\hfil \SBI #2}
 \setbox4=\hbox to \wd0{\hfil \SBI #3}
 \setbox6=\hbox to \wd0{\hss
             \vbox{\hsize=\wd0 \parskip=0pt \baselineskip=10 true pt
                   \copy0 \break%
                   \copy2 \break% 
                   \copy4 \break}}
 \dimen0=\ht6   \advance\dimen0 by \vsize \advance\dimen0 by 8 true pt
                \advance\dimen0 by -\pagetotal
	        \advance\dimen0 by \IMSmarkvadjust
 \dimen2=\hsize \advance\dimen2 by .25 true in
	        \advance\dimen2 by \IMSmarkhadjust

%
%   Check for publication info
%
%  \newread\jref
  \openin2=publishd.tex
  \ifeof2\setbox0=\hbox to 0pt{}
  \else 
     \setbox0=\hbox to 3.1 true in{
                \vbox to \ht6{\hsize=3 true in \parskip=0pt  \noindent  
                {\SBI \IMSpubltext}\hfil\break
                \input publishd.tex 
                \vfill}}
  \fi
  \closein2
  \ht0=0pt \dp0=0pt
 \ht6=0pt \dp6=0pt
 \setbox8=\vbox to \dimen0{\vfill \hbox to \dimen2{\copy0 \hss \copy6}}
 \ht8=0pt \dp8=0pt \wd8=0pt
 \copy8
 \message{*** Stony Brook IMS Preprint #1, #2. #3 ***}
}
 
\def\IMSmarkvadjust{-30pt}
\SBIMSMark{2002/03}{October 2002}{}

\begin{abstract}
In this work, under a mild assumption, we give the
classification of the complete polynomial
vector fields in two variables up to algebraic automorphisms
of $\C^2$. The general problem is also reduced to the study of
the combinatorics of certain resolutions of singularities.
Whereas we deal with $\C$-complete vector fields,
our results also apply to $\R$-complete ones thanks to a
theorem of Forstneric \cite{forst}.
\end{abstract}

\bigskip

\noindent \hspace{0.9cm} {\small Key-words: line at infinity  - vector fields - singular
foliations}

\bigskip

\noindent \hspace{0.9cm} {\small AMS-Classification  34A20, 32M25, 32J15}

\section{Introduction}

\hspace{0.4cm}
Recall that a {\it holomorphic flow} on $\C^2$ is a holomorphic mapping
$
\Phi \, : \; \C \times \C^2 \longrightarrow \C^2
$
satisfying the two conditions below:

\noindent $\bullet$ $\Phi (0, p) = p$ for every $p \in \C^2$;

\noindent $\bullet$ $\Phi (T_1 + T_2 ,p ) = \Phi (T_1 ,
\Phi (T_2, p))$.

A holomorphic flow $\Phi$ on $\C^2$ induces a {\it holomorphic
vector field} $X$ on $\C^2$ by the equation
$$
X (p) = \left. \frac{d \Phi (T,p)}{dT} \right|_{T=0} \, .
$$
Conversely a holomorphic vector field $X$ on $\C^2$
is said to be {\it complete} if it is associated to a 
holomorphic flow $\Phi$. Since every polynomial vector
field of degree~$1$ is complete, we assume that $X$
has degree~$2$ or greater. A polynomial vector field $X$ can be
considered as a meromorphic vector field on $\pp$ therefore
inducing a singular holomorphic foliation $\fol_X$
on $\pp$. The singularities of $\fol_X$ lying in the ``line at
infinity'' $\Delta$ will be denoted by $p_1 ,\ldots ,p_k$.
A singularity $p_i$ as above is called {\it dicritical} if
there are infinitely many analytic curves invariant
by $\fol_X$ and passing through $p_i$.
The first result of this paper
is the following:

\bigskip

\noindent {\bf Theorem A} {\sl Let $X$ be a complete polynomial
vector field on $\C^2$ with degree~$2$ or greater and let
$\fol_X$ be the singular foliation induced by $X$ on $\pp$. Assume
that $\fol_X$ has a dicritical singularity in $\Delta$.
Then $X$ is conjugate by a polynomial
automorphism to
one of the following vector fields:
\begin{enumerate}

\item $P(y) x^{\epsilon} \partial /\partial x$, where $P(y)$ is a polynomial
in $y$ and $\epsilon=0,\, 1$.

\item $x^ny^m(mx \parlx - ny \parly)$, where ${\rm g.c.d} (m,n) =1$ and $m,n \in \N$;

\end{enumerate}
}

In view of Theorem~A, we just need to consider vector fields
all of whose singularities belonging to $\Delta$ are not dicritical.
Again let $X$ be such a vector field and let $\fol_X$ be its associated
foliation.
Consider a singularity $p_i$ of $\fol_X$ in the line at infinity
$\Delta$ and a vector field $\widetilde{X}$ obtained through
a finite sequence of blowing-ups of $X$ beginning at $p_i$. Denote
by $\cale$ the corresponding exceptional divisor and by $D_i$, $i=1,
\ldots ,l$,
its irreducible components which are all rational curves. We say
that $X$ has adapted poles at $p_i$ if, for every sequence
of blow-ups as above and every irreducible component $D_i$ of
the corresponsding exceptional divisor
$\cale$, either $X$ vanishes identically on $D_i$ or $D_i$ consists of pole
of $\widetilde{X}$ (in other words $\widetilde{X}$ is not regular
on $D_i$). Clearly the great majority of polynomial vector fields
have adapted poles at its singularities at infinity. We then have:

\bigskip

\noindent {\bf Theorem B} {\sl Let $X$ be a complete polynomial
vector field on $\C^2$ and denote by $p_i$, $i=1, \ldots ,k$
the singularities of the associated foliation $\fol_X$ belonging
to the line at infinity $\Delta$. Suppose that $X$ has adapted
poles at each $p_i$. Then $\fol_X$ possesses a dicritical singularity
in the line at infinity $\Delta$.}

There is a good amount of literature devoted to complete
vector fields on $\C^2$, in particular
our results are complementary to the recent results obtained
by Cerveau and Scardua in \cite{cesc}. Note however that
the points of view adopted
in both papers are almost disjoint. For more information
on complete polynomial vector fields the reader can consult
the references at the end as well as references
in \cite{cesc}.

After Theorems~A and~B, in order to classify all complete
polinomial vector fields on $\C^2$ we just have to consider
vector fields which do not have adapted poles at one of the
singularities $p_1, \ldots , p_k$. In particular we can always
assume that none of these singularities is dicritical.

Let us close this Introduction by giving a brief description of the structure
of the paper. First we observe that the
method developed here may be pushed forward to deal
with vector fields which do not have adapted poles.
Indeed the assumption  that $X$ has
adapted poles in $\Delta$ is used only in Section~6.
More precisely, from Section~2 to Section~5, the classification
of complete polynomial vector fields is reduced to a problem of
understanding the possible configurations of rational curves
arising from blow-ups of the singularities of $\fol_X$ in the
line at infinity $\Delta$. The role of our main assumption is
to make the ``combinatorics'' of these configurations simpler
so as to allow for the complete description given in Section~6.
It is reasonable to expect that a more detailed study of these
configurations will lead to the general classification
of complete polynomial vector fields.

Another feature of our method is its local nature. Indeed most
of our results are local and therefore have potential to be
applied in other situations (especially to other Stein surfaces).
We mention in particular the results of Sections~3
and~5 (cf. Theorem~(\ref{selano}), Proposition~(\ref{prop4.2})).
Also the local vector fields $Z_{1,11}, \, Z_{0,12}, \,
Z_{1,00}$ introduced in Proposition~(\ref{prop4.2}) might
have additional interest.

This paper is also a sequence of \cite{re3} where it was observed,
in particular, that the problem of understanding complete
polynomial vector fields on $\C^2$ can be unified with
classical problems in Complex Geometry through the notion
of {\it meromorphic semi-complete vector fields}.
The method employed in the proof of our main result relies
heavily on this connection. Indeed an important part of the proof
the preceding theorems is a discussion of semi-complete singularities
of meromorphic vector fields. The study of these singularities was
initiated in \cite{re3} but the present discussion is based on
a different setting.

\noindent {\bf Acknowledgements}: I am grateful to D. Cerveau and
B. Scardua who raised my interest in this problem by sending me
their preprint \cite{cesc}.

\section{Basic notions and results}

\hspace{0.4cm}
The local orbits of a polynomial vector field $X$
induce a singular holomorphic foliation $\fol_X$ on $\C^2$.
Besides, considering $\pp$ as the natural compactification
of $\C^2$ obtained by adding the line at infinity $\Delta$,
the foliation $\fol_X$ extend to a holomorphic foliation,
still denoted by $\fol_X$, on the whole of $\pp$. 
This extension may or may not leave the line at infinity $\Delta$
invariant.
On the other hand, the vector field $X$ possesses a
{\it meromorphic} extension to $\pp$, also denoted by $X$,
whose {\it pole divisor} coincides with $\Delta$. Note that
the meromorphic extension of $X$ to $\pp$ happens to be holomorphic
if and only if the degree of $X$ is $1$ or if $X$ has degree
$2$ and the line at infinity $\Delta$ is not invariant
by $\fol_X$ (for further details cf. below).
Let us
make these notions more precise.

Recall that a
{\it meromorphic}\, vector field $Y$ on a neighborhood $U$
of the origin $(0, \ldots , 0) \in \C^n$ is by definition
a vector field of the form
$$
Y = F_1 \frac{\partial}{\partial z_1} + \cdots +
F_n \frac{\partial}{\partial z_n} \, ,
$$
where the $F_i$'s are meromorphic functions on $U$
(i.e. $F_i = f_i /g_i$ with $f_i ,g_i$ holomorphic on
$U$). Note that $Y$ may not be defined on the whole $U$
even though we consider $\infty$ as a value since
$F_i$ may have indeterminacy points. We denote by
$D_Y$ the union of the sets $\{ g_i =0 \}$. Of course
$D_Y$ is a divisor consisting of poles and
indeterminacy points of $Y$ which is called
the {\it pole divisor} of $Y$

\begin{defnc}
The meromorphic vector field $Y$ is said to be semi-complete on $U$
if and only if there exists a meromorphic map
$\Phi_{sg} : \Omega \subseteq \C \times U \rightarrow
U$, where $\Omega$ is an open set of
$\C \times U$, satisfying the conditions below.
\begin{enumerate}

\item $$
\left. \frac{d \Phi_{sg} (T, x)}{dT} \right|_{T=0}
= Y(x) \; \, \mbox{for all $x \in U \setminus D_x$;}
$$

\item $\Phi_{sg} (T_1 + T_2, x) = \Phi_{sg} (T_1,
\Phi_{sg} (T_2 ,x))$ provided that both
sides are defined;

\item If $(T_i ,x)$ is a sequence of points in $\Omega$
converging to a point $(\hat{T} ,x)$ in the boundary
of $\Omega$, then $\Phi_{sg} (T_i ,x)$ converges to
the boundary of $U \setminus D_Y$ in the sense that
the sequence leaves every compact subset of
$U \setminus D_Y$.
\end{enumerate}
\end{defnc}

The map $\fii$ is called the meromorphic
semi-global flow associated to $Y$ (or induced by $Y$).

Assume we are given a meromorphic vector field $Y$ defined on a neighborhood
of $(0,0) \in \C^2$. It is easy to see that $Y$ has the form $Y = f Z/g$
where $f,g$ are holomorphic functions and $Z$ is a holomorphic
vector field having at most an isolated singularity at the origin.
Naturally we suppose that $f,g$ do not have a (non-trivial)
common factor, so that
$f,g$ and $Z$ are unique (up to a trivial, i.e. inversible, factor).
Next, let $\fol$ denote the local singular foliation
defined by the orbits of $Z$. We call $\fol$ the foliation
associated to $Y$ and note that either $\fol$ is regular at
the origin or the origin is an isolated singularity of $\fol$.
An analytic curve $\cals$ passing through the origin and
invariant by $\fol$ is said to be a {\it separatrix} of $\fol$
(or of $Y, Z$).

%The theorem below plays the main role in the proof of the
%theorem stated in the Introduction. Furthermore this
%theorem is likely to have further applications.

%\begin{teo}
%\label{main}
%Let $Y$ be a meromorphic vector field defined and semi-complete
%on a neighborhood of the origin in $\C^2$ and consider
%the associated foliation $\fol$. Set $Y = f Z/g$
%as above and suppose that $\{ g=0\}$ is a smooth (irreducible)
%separatrix of $\fol$......................
%Then $Y$ admits one of the following
%holomorphic normal forms:
%\end{teo}

The rest of this section is devoted to establishing some
preliminary results concerning both theorems in the
Introduction. Particular attention will be paid to meromorphic
semi-complete vector fields which appear when we restrict $X$
to a neighborhood of the line at infinity $\Delta$. As it will
be seen, a large amount of information on $X$ arises from a detailed
study of this restriction.
To begin with, let us
recall the notion of {\it time-form} $dT$ of a meromorphic
vector field and the basic lemma about integrals
of $dT$ over curves. For other general facts
about meromorphic semi-complete vector fields
the reader is referred to Section~2 of \cite{re3}.

Let $Y$ be a meromorphic vector field defined on an open
set $U$ and let $\fol$ denote its associated foliation. The regular
leaves of $\fol$ (after excluding possible punctures corresponding
to {\it zeros} or poles of $Y$) are naturally equipped with
a {\it foliated} holomorphic $1$-form $dT$ defined by imposing
$dT . Y = 1$. As a piece of terminology, whenever the $1$-form $dT$
is involved, the expression ``regular leaf of $\fol$'' should be
understood as
a regular leaf $L$ (in the sense of the foliation $\fol$) from
which the intersection with the set of {\it zeros} or poles
of $Y$ was deleted. Hence the restriction of $dT$ to a regular
leaf $L$ is, by this convention, always holomorphic.
Note also that $dT$ is ``foliated'' in the sense
that it is defined only on the tangent spaces of the leaves.
We call $dT$ the {\it time-form} associated to (or induced by) 
$Y$. Lemma~(\ref{timeform}) below is the most fundamental result
about semi-complete vector fields (cf. \cite{re1}, \cite{re3}).

\begin{lema}
\label{timeform}
Let $Y$, $U$, $\fol$ and $dT$ be as above. Consider a regular leaf
$L$ of $\fol$ and an embedded (open) curve $c: [0,1]
\rightarrow L$. If $Y$ is semi-complete on $U$ then the
integral of $dT$ over $c$
does not vanish.\qed
\end{lema}

Let us now go back to a complete polynomial vector field $X$
on $\C^2$ whose degree is $d \in \N$. Set
\begin{equation}
X = X_0 + X_1 + \cdots + X_d \label{equa1}
\end{equation}
where $X_i$, $i =1 ,\ldots ,d$, stands for the homogeneous
component of degree $i$ of $X$. With this notation
the vector fields whose associated foliations $\fol_X$ do not
leave the line at infinity $\Delta$ invariant
admit an elementary characterization
namely: $\fol_X$ does not leave $\Delta$ invariant if and only if $X_d$ has
the form $F(x,y) (x \parlx \, + \, y \parly)$, where
$F$ is a homogeneous polynomial of degree $d-1$.
Furthermore, viewing $X$
as a meromorphic vector field on $\pp$, a direct inspection shows
that the
order of the pole divisor $\Delta$ is $d-1$
provided that $\Delta$ is invariant under $\fol_X$. If $\Delta$
is not invariant under $\fol_X$ then this order is $d-2$.
On the other hand, given a point $p \in \Delta$
and a neighborhood $U \subset \pp$ of $p$,
it is clear that $X$ defines a meromorphic semi-complete vector
field on $U$.

Our first lemma shows that we can suppose that
the line at infinity is invariant
by the associated foliation $\fol_X$. Whereas the proof
is elementary, we give a detailed account since some basic
ideas will often be used later on.

\begin{lema}
\label{lema2.1}
Consider a complete polynomial vector field $X$ on
$\C^2$ and denote by $\fol_X$ the foliation induced by
$X$ on $\pp$. Assume that the line at infinity $\Delta$ is not
invariant under $\fol_X$. Then the degree
of $X$ is at most $1$.
\end{lema}

\noindent {\it Proof}\,: First we set $X = F Z$ where
$F$ is a polynomial of degree $0 \leq n \leq d$ and
$Z$ is a polynomial vector field of degree $d-n$ and
isolated {\it zeros}. In other words, we have
$Z = P \parlx + Q \parly$ where $P, Q$ for polynomials $P,Q$ without
non-trivial common factors.

First we
suppose for a contradiction that $d$ is strictly greater
than $2$. In view of the preceding discussion, the line
at infinity $\Delta$ is the polar divisor of $X$ and has order
$d-2 \geq 1$. Let ${\mathcal C} \subset \pp$ be the algebraic
curve induced in affine coordinates by $F=0$. Finally
consider a ``generic'' point $p \in \Delta$ and a neighborhood $U$
of $p$ such that $U \cap {\mathcal C} = \emptyset$.

Let $\fol_X$ be the singular foliation induced by $X$
on $\pp$ and notice that $p$ is a regular point of $\fol_X$.
Besides the leaf $L$ containing $p$ is transverse at $p$
to $\Delta$. Thus we can introduce coordinates $(u,v)$ on $U$,
identifying $p$ with $(0,0) \in \C^2$, in which
$X$ becomes
$$
X(u,v) = u^{2-d}. f. \frac{\partial}{\partial u} 
$$
where $f$ is a holomorphic function such that
$f(0,0) \neq 0$ and $\{ u = 0 \} \subset \Delta$
(here we use the fact that $U \cap {\mathcal C} =
\emptyset$).
The axis $\{ v=0 \}$ is obviously invariant under
$\fol$ and the time-form $dT_{\{v = 0\}}$ induced
on $\{ v =0 \}$ by $X$ is given by
$dT_{\{v = 0\}} = u^{d-2} du / f(u,0)$. Since
$d \geq 3$ and $f(0,0) \neq 0$, it easily follows
the existence of an embedded curve $c:[0,1]
\rightarrow \{ v=0 \} \setminus (0,0)$ on which
the integral of $dT_{\{v = 0 \}}$ (cf. Remark~(\ref{obsidiota})).
The resulting
contradiction shows that $d \leq 2$.

It only remains to check the case $d=2$. Modulo
a linear change of coordinates, we have
$X= X_0 + X_1 + x( x \parlx \, + \, y \parly)$.
The above calculation shows that the natural extension
of $X$ to $\pp$ is, in fact, holomorphic. Therefore
$X$ is complete on $\pp$ since $\pp$ is compact.
Besides a generic point $p$
of $\Delta$ is regular for $X$ and has its
(local) orbit transverse to $\Delta$.
It follows that points in the orbit of $p$
reaches $\Delta$ in finite time. Thus the restriction
of $X$ to the affine $\C^2$ cannot be complete
as its flow goes off to infinity in finite time.\qed

In view of Lemma~(\ref{lema2.1}) all complete polynomial
vector fields considered from now on will be such that the
associated foliation $\fol_X$ leaves the line at infinity
$\Delta$ invariant. In particular, in the sequel,
the extension of $X$ to $\pp$ is strictly meromorphic.
Furthermore the pole divisor is constituted by the
line at infinity $\Delta$ and has order $d-1$ (where
$d$ is the degree of $X$).

\begin{lema}
\label{lema2.2}
Let $X$ be as above and let $X_d$ be its top-degree homogeneous
component (as in (\ref{equa1})). Then $X_d$ is semi-complete
on the entire $\C^2$.
\end{lema}

\noindent {\it Proof}\,: For each integer $k \in \N$, we consider
the homothety $\Lambda_k (x,y) = (k x , k y)$ of $\C^2$. The vector
fields defined by $\Lambda_k^{\ast} X$ are obviously
complete, and therefore semi-complete, on $\C^2$. Next we set
$X^k = k^{1-d} \Lambda_k^{\ast} X$. Since $X^k$ and
$\Lambda_k^{\ast} X$ differ by a multiplicative constant,
it is clear that $X^k$ is complete on $\C^2$. Finally, when
$k$ goes to infinity, $X^k$ converges uniformly on $\C^2$
towards $X_d$. Since the space of semi-complete vector fields
is closed under uniform convergence (cf. \cite{ghre}), it results
that $X_d$ is semi-complete on $\C^2$. The lemma is proved.\qed

The lemma above has an immediate application. In fact, a
homogeneous polynomial vector field has a $1$-parameter
group of symmetries consisting of homotheties. Hence these
vector fields can essentially be integrated. In other
words, it is possible to describe all homogeneous polynomial
vector fields which are semi-complete on $\C^2$. This
classification was carried out in \cite{ghre} and, combined
with Lemma~(\ref{lema2.2}), yields:

\begin{coro}
\label{lema2.3}
Let $X$ and $X_d$ be as in the statement of Lemma~(\ref{lema2.2}).
Then, up to a linear change of coordinates of $\C^2$, $X_d$
has one of the following normal forms:
\begin{enumerate}

\item $X_d = y^a f(x,y) \parlx$ where $f$ has degree strictly less
than~$3$ and $a \in \N$.

\item $X_d = x(x \parlx + ny \parly)$, $n \in 
N$, $n\neq 1$.

\item $X_d = x^i y^j (mx \parlx - n y \parly)$ where
$m,n \in \N^{\ast}$ and $ni-mj = -1, \, 0 \, 1$. 
We also may have $X_d = (xy)^a (x \parlx - y \parly)$,
$a \in \N$.

\item $X_d = x^2 \parlx - y(nx-(n+1)y) \parly$, $n \in \N$.

\item $X_d = [xy(x-y)]^a [x(x-2y) \parlx + y(y-2x) \parly]$,
$a \in \N$.

\item $X_d = [xy(x-y)^2]^a [x(x-3y) \parlx + y(y-3x) \parly]$,
$a \in \N$.

\item $X_d = [xy^2(x-y)^3]^a [x(2x-5y) \parlx + y(y-4x) \parly]$,
$a \in \N$.\qed
\end{enumerate}
\end{coro}
As an application of
Corollary~(\ref{lema2.3}), we
shall prove Lemma~(\ref{lema2.4}) below. This lemma estimates
the number of singularities
that the foliation $\fol_X$ induced by $X$ on $\pp$ may have
in the line at infinity.
This estimate will be useful in Section~7. Also recall that
the line at infinity is supposed to be invariant by $\fol_X$
(cf. Lemma~(\ref{lema2.1})).

\begin{lema}
\label{lema2.4}
Let $X$ be a complete polynomial vector field on $\C^2$
and denote by $\fol_X$ the foliation induced by $X$ on $\pp$.
Then the line at infinity contains
at most $3$ singularities of $\fol_X$.
\end{lema}

\noindent {\it Proof}\,: We consider the change of coordinates
$u=1/x$ and $v = y/x$ and note that in the coordinates $(u,v)$ the line
at infinity $\Delta$ is represented by $\{ u=0\}$.

First we set $X_d = F. Y_d$ where $F$ is a polynomial
of degree $k$ and $Y_d$ is a polynomial vector field of
degree $d-k$.  Next let us consider the algebraic curve
${\mathcal C} \subset \pp$ induced on $\pp$ by the affine
equation $F=0$. We also consider the foliation $\fol_d$ induced
on $\pp$ by $Y_d$. Finally denote by $\Delta \cap {\rm Sing} (\fol_X)$
(resp. $\Delta \cap {\rm Sing} (\fol_d)$) the set of singularities
of $\fol_X$ (resp. $\fol_d$) belonging to $\Delta$.
An elementary calculation with the coordinates
$(u,v)$ shows that
$$
\Delta \cap {\rm Sing} (\fol_X) \subseteq (\Delta \cap {\rm Sing} (\fol_d))
\cup ({\mathcal C} \cap \Delta) \, .
$$
Now a direct inspection in the list of Corollary~(\ref{lema2.3})
implies that the set $(\Delta \cap {\rm Sing} (\fol_d))
\cup ({\mathcal C} \cap \Delta)$ consists of at most~$3$
points. The proof of the lemma is over.\qed

\section{Simple semi-complete singularities}

\hspace{0.4cm}
In this section we shall begin the study of a certain
class of semi-complete singularities. The results obtained
here will largely be used in the remaining
sections. In the sequel $Y$ stands for
a meromorphic vector field defined on a neighborhood of $(0,0)
\in \C^2$. We always set $Y = fZ/g$ where $f,g$ are holomorphic
functions without common factors and $Z$ is a holomorphic
vector field for which the origin is either a regular point
or an isolated singularity. Also $\fol$ will stand for the
foliation associated to $Y$ (or to $Z$). We point out that the
decomposition $Y = fZ/g$ is unique up to an inversible factor.

If $Z$ is singular at $(0,0)$, we can consider its eigenvalues
at this singularity. Three cases can occur:

\noindent {\bf a-} Both eigenvalues, $\lambda_1 , \lambda_2$
of $Z$ vanish at $(0,0)$.

\noindent {\bf b-} Exactly one eigenvalue, $\lambda_2$, vanishes
at $(0,0)$.

\noindent {\bf c-} None of the eigenvalues $\lambda_1 , \lambda_2$
vanishes at $(0,0)$

\noindent Whereas $Z$ is defined up to an inversible factor, all the cases
{\bf a}, {\bf b} and {\bf c} are well-defined. In the case {\bf c},
the precise values of $\lambda_1 , \lambda_2$ are not well-defined
but so is their ratio $\lambda_1 / \lambda_2$. Following a usual
abuse of notation, in the case {\bf c}, we shall say that 
the foliation $\fol$ associated to $Z$
has eigenvalues $\lambda_1 , \lambda_2$ different from {\it zero}.
In other words, given a singular holomorphic foliation $\fol$,
we say that $\fol$ has eigenvalues $\lambda_1 ,\lambda_2$
if there exists $Z$ as before having $\lambda_1 ,\lambda_2$
as its eigenvalues at $(0,0)$.
The reader will easily check that all the relevant notions
discussed depend only on the ratio $\lambda_1 / \lambda_2$.
A singularity is said to be {\it simple} if it has at least one
eigenvalue different from {\it zero}. A simple singularity possessing
exactly one eigenvalue different from {\it zero} is called a
{\it saddle-node}.

More generally the {\it order} of $\fol$ at $(0,0)$ is
defined as the order of $Z$ at $(0,0)$, namely it is
the degree of the first non-vanishing homogeneous component
of the Taylor series of $Z$ based at $(0,0)$. It is obvious
that the order of $\fol$ does not depend on the vector
field with isolated singularities $Z$ chosen.

\begin{obs}
\label{obsidiota}
{\rm A useful fact is the non-existence, in dimension~$1$,
of {\it strictly
meromorphic} semi-complete vector fields. In other words,
if $Y = f(x) \partial /\partial x$ with $f$ meromorphic,
then $Y$ is not semi-complete on a neighborhood of $0 \in \C$.
In fact, fixed a neighborhood $U$ of $0 \in \C$,
we have $f(x) = x^{-n}. h(x)$ where $n \geq 1$,
$h$ is holomorphic and $h (0) \neq 0$. Thus the corresponding
time-form is $dT = x^n dx /h$. It easily follows the existence
of an embedded curve $c : [0,1] \rightarrow U \setminus \{ (0,0) \}$
on which the integral the integral of $dT$ vanishes.
Similarly we can also prove that $Y$ is not semi-complete
provided that $0 \in \C$ is an essential singularity of $f$.

Summarizing the preceding
discussion, the fact that $Y$ is semi-complete implies that
$f$ is holomorphic at $0 \in \C$. Consider now that $f$
is holomorphic but $f(0) = f'(0) = f''(0) =0$.
An elementary estimate (cf. \cite{re1}) shows that in this case
$Y$ is not semi-complete. Finally when $Y$ is semi-complete
but $f(0) = f'(0) =0$ it is easy to see that $Y$ is conjugate to
$x^2 \parlx$ (cf. \cite{ghre}). These elementary results give
a complete description of semi-complete singularities in dimension~$1$.}
\end{obs}

Let us say that $P = P_{\alpha} /P_{\beta}$ is a {\it homogeneous
rational function} if $P_{\alpha}$ and $P_{\beta}$ are homogeneous
polynomials (possibly with different degrees).
The next lemma is borrowed from \cite{re3}

\begin{lema}
\label{le3.1}
Consider the linear vector field $Z = x \partial
/\partial x + \lambda y \partial /\partial y$.
Suppose that $P = P_{\alpha} /P_{\beta}$ is a (non-constant)
homogeneous rational function
and that $\lambda \not\in \R_{+}$. Suppose also that
$P Z$ is semi-complete. Then one has

\noindent 1. $\lambda$ is rational, i.e. $\lambda =-n/m$
for appropriate coprime positive integers $m,n$.

\noindent 2. $P$ has one of the forms below:

\noindent {\bf 2i} $P = x^c y^d$ where $mc-nd=0$ or $\pm 1$.

\noindent {\bf 2ii} If $\lambda =-1$, then $P$ is
$x^c y^d$ ($mc-nd=0$ or $\pm 1$) or
$P = (x-y)(xy)^a$ for $a \in \Z$.\qed
\end{lema}

\begin{obs}
\label{ipc}
{\rm Consider a holomorphic vector field of the form
$$
x^a y^b h(x,y) [ x(1 + {\rm h.o.t.}) \parlx \; - \;
y(1 + {\rm h.o.t.}) \parly \, ,
$$
where $a,b \in \Z$ and $h$ is holomorphic with $h(0,0) =0$.
Of course we suppose that $h$ is not divisible by
$x,y$. Next we assume that $X$ is semi-complete on
a neighborhood of $(0,0)$. Denote by $h^k$ the
homogeneous component of the first non-trivial
jet of $h$ at $(0,0)$. The same argument employed in the
proof of Lemma~(\ref{lema2.2}), modulo replacing $k$ by
$1/k$, shows that the vector field
$x^a y^b h^k (x \parlx - y \parly)$ is semi-complete.
From the preceding lemma it then follows that $h^k = x-y$
and $a=b$. However a much stronger conclusion holds: $X$
admits the normal form
$$
(xy)^a (x-y) g(x \parlx - y \parly) \, ,
$$
where $g$ is a holomorphic function satisfying $g(0,0) \neq 0$.
In fact, in order to deduce the normal form above, we just
need to check that $x(1 + {\rm h.o.t.}) \parlx \; - \;
y(1 + {\rm h.o.t.}) \parly$ is linearizable. After \cite{mamo},
this amounts to prove that the local holonomy of their
separatrizes is the identity. However the integral of time-form
on a curve $c$ projecting onto a loop around $0$ in $\{ y=0\}$
is cearly equal to {\it zero}. Since $X$ is semi-complete, such curve must
be closed which means that the holonomy in question is trivial.}
\end{obs}

Of course the next step is to discuss the case $\lambda >0$.
However, at this point, we do not want to consider only linear
vector fields. This discussion will naturally lead us to
consider singularities having an infinite number of separatrizes.
Recall that a singularity of a holomorphic foliation $\fol$ is said to be
{\it dicritical} if $\fol$ possesses infinitely many
separatrizes at $p$. Sometimes we also say that a vector field
$Y$ defined on a neighborhood of $(0,0) \in \C^2$ is dicritical
to say that $(0,0)$ is a dicritical singularity of the foliation
associated to $Y$. Let us begin with the following:

\begin{lema}
\label{le3.2}
Consider a semi-complete meromorphic vector field $Y$ defined on a
neighborhood of $(0,0) \in \C^2$ and having the form
$$
Y = \frac{f}{g} Z
$$
where $f,g$ are holomorphic functions with $f(0,0)g(0,0)=0$
and $Z$ is a holomorphic vector field with an isolated
singularity at $(0,0) \in \C^2$ whose eigenvalues are
$1$ and $\lambda$.
Assume that $\lambda >0$ but neither $\lambda$ nor $1/\lambda$
belongs to $\N$.
Then $Z,Y$ are dicritical vector fields.
\end{lema}

\noindent {\it Proof}\,: Note that
Poincar\'e's linearization theorem
\cite{ar} ensures that $Z$ is linearizable. Therefore, in appropriate
coordinates, we have $Y =f(x \parlx + \lambda y \parly) / g$.
If $\lambda$ is rational equal to $n/m$, then $Z$ admits
the meromorphic first integral $x^n y^{-m}$ and therefore
admits an infinite number of separatrizes.

In order to prove the lemma is now sufficient to check that
$\lambda$ is rational provided that $Y$ is semi-complete.
Let $P_{\alpha}$ (resp. $P_{\beta}$) be the
first non-vanishing homogeneous component of the Taylor
series of $f$ (resp. $g$) at $(0,0) \in \C^2$. The same
argument carried out in the proof of Lemma~(\ref{lema2.2}),
modulo replacing $k$ by $1/k$, shows that the vector
field $Y^{\rm ho} = P_{\alpha} (x \parlx + \lambda y \parly) /P_{\beta}$
is semi-complete on $\C^2$. We are going to see that this
implies that $\lambda$ is rational.

Suppose for a contradiction that $\lambda$ is not rational.
In this case the only separatrizes of $Y^{\rm ho}$ are
the axes $\{ x=0 \}$, $\{ y=0\}$. Since in dimension~$1$ there is no
meromorphic semi-complete vector field, it
follows that the {\it zero set} of $P_{\beta}$
has to be invariant under the foliation $\fol_{Y^{\rm ho}}$ associated
to $Y^{\rm ho}$. Thus $P_{\beta}$
must have the form $x^a y^b$ for some $a,b \in \N$. Therefore
we can write $P$ as $x^cy^d Q(x,y)$ where $Q$ is a homogeneous
polynomial.

Observe that the orbit $L$ of $Y^{\rm ho}$ (or $\fol_{Y^{\rm ho}}$)
passing through
the point $(x_1 ,y_1)$ ($x_1 y_1 \neq 0$) is parametrized
by $\textsf{A} : \, T \mapsto (x_1 e^T , y_1 e^{\lambda T})$.
The restriction to $L$ of the vector field $P Z$
is given in the coordinate $T$ by
$P(x_1 e^T , y_1e^{\lambda T}) \partial /\partial T$.
Because $\lambda$ is not rational, the parametrization $\textsf{A}$
is an one-to-one map from $\C$ to $L$. It
results that the one-dimensional vector field
$x_1^cy_1^d e^{(c+\lambda d)T} Q(x_1e^T ,
y_1^{\lambda T}) \partial /\partial T$ is semi-complete
on the entire $\C$. On the other hand the function $T \mapsto
e^{(c+\lambda d)T} Q(x_1e^T , y_1^{\lambda T})$
is defined on the whole of $\C$. Since $\lambda$
is not rational and $Q$ is a polynomial, we conclude
that this function has an essential singularity at 
infinity. This
contradicts the fact that this function
corresponds to a semi-complete vector field
(cf. Remark~(\ref{obsidiota})). The lemma is proved.\qed

Let us now consider the case $\lambda \in \N$ since the case
$1/\lambda \in \N$ is analogous.
Thus we denote by
$1$ and $n \in \N^{\ast}$ the eigenvalues of $Z$ at $(0,0)$.
Such $Z$ is either linearizable
or conjugate to its Poincar\'e's-Dulac normal form \cite{ar}
\begin{equation}
(nx + y^n) \parlx \; + \; y \parly \, . \label{awui}
\end{equation}
When $Z$ is linearizable, it has infinitely
many separatrizes. Thus we are, in fact, 
interested in the case in which $Z$ is
conjugate to the Poincar\'e-Dulac normal form~(\ref{awui}).
In particular $\{ y =0 \}$ is the unique separatrix
of $Z$ (or $Y$).

\begin{lema}
\label{le3.3}
Let $Y$ be $Y = fZ/g$ where $Z$ is a vector field as
in~(\ref{awui}) and $f,g$ are holomorphic functions
satisfying $f(0,0) g (0,0) =0$. Assume that $Y$ is semi-complete
and that the regular orbits of $Y$ can contain at most one
singular point whose order is necessarily~$1$. Then $n=1$.
Furthermore, up to an inversible factor, $g(x,y) = y$ and $\{ f=
0\}$ defines a smooth analytic curve which is not tangent
to $\{ y=0\}$.
\end{lema}

\noindent {\it Proof}\,: Since the divisor of poles of $Y$ is contained
in the union of the separatrizes, it follows that $Y$ has the form
$$
Y = y^{k} F(x,y) [(nx + y^n) \parlx \; + \; y \parly] \, ,
$$
where $k \in \Z$ and $F$ is a holomorphic function which is
not divisible by $y$. Clearly $k < 0$, otherwise the
first homogeneous component of $Y$ would not be semi-complete
(cf. Corollary~(\ref{lema2.3}) and Remark~(\ref{obsidiota})).

Let us first deal with the case $n=1$. Note that we are going
to strongly use Theorem~(\ref{selano})
which is the next result to be proved. This theorem is concerned
with the so-called {\it saddle-node} singularities
which are those having exactly one eigenvalue different from
{\it zero}.
Blowing-up the vector field $Y$ we obtain a new vector field
$\wiX$ defined and semi-complete on a neighborhood
of the exceptional divisor $\pi^{-1} (0)$ (where $\pi$ stands
for the blow-up map). Denote by $\wifol$ the foliation associated
to $\wiX$ and note that $\wifol$ has a unique singularity
$p \in \pi^{-1} (0)$. More precisely, $\wiX$ on a neighborhood of $p$
is given in standard coordinates $(x,t)$ ($\pi (x,t) = (x, tx)$),
by
$$
H (x,t) [(x(1+t) \parlx \; - \; t^2 \partial /\partial t] \, ,
$$
where $H$ is meromorphic function. Because of Theorem~(\ref{selano}),
we know that the restriction of $\wiX$ to the exceptional
divisor $\{ x=0\}$ has to be regular i.e. $H$ is not divisible
by $x$ or $x^{-1}$. This implies that the order of $F$ at
$(0,0) \in \C^2$ is $k$ and, in particular, $F(0,0) =0$.

On the other hand, $H$ has the form $t^k h (x,t)$ where
$h$ is holomorphic on a neighborhood of $x=t=0$ and not divisible
by $t$ or $t^{-1}$. Again Theorem~(\ref{selano}) shows that
$h$ has to be an inversible factor, i.e. $h(0,0) \neq 0$. In
other words, the proper transform of the (non-trivial) analytic
curve $F=0$ intersects $\pi^{-1} (0)$ at points different from
$x=t=0$.

The restriction of $\wiX$ to $\pi^{-1} (0)$ is a holomorphic
vector field which has a singularity at $\{x=t=0\}$ whose order
is $2-k$. In particular $k \in \{ 0, 1,2\}$. The other singularities
correspond to the intesection of $\pi^{-1} (0)$ with the
proper transform of $\{ F= 0\}$. The statement follows
since the regular orbits of $X$ can contain only one singular 
point whose order is~$1$.\qed

In the rest of this section we briefly discuss the case
of singularities as in {\bf b}, that is, those singularities
having exactly one eigenvalue different from {\it zero}.
As mentioned
they are called {\it saddle-nodes}
and were classified in \cite{mara}. A consequence
of this classification is the existence of
a large moduli space. The subclass of saddle-nodes consisting of those
associated to semi-complete holomorphic
vector fields was characterized in \cite{re4}. In the sequel we
summarize and adapt these results to meromorphic semi-complete
vector fields.

To begin with, let $\omega$ be a singular holomorphic
$1$-form defining a saddle-node $\fol$. According to
Dulac \cite{dulac}, $\omega$ admits the normal form
$$
\omega (x,y) = [x(1+ \lambda y^p) + yR(x,y)] \, dy \; - \;
y^{p+1} \, dx \; ,
$$
where $\lambda \in \C$, $p \in \N^{\ast}$ and $R(x,0) =
o (\mid x \mid^p)$. In particular $\fol$ has a (smooth)
separatrix given in the above coordinates by $\{ y=0\}$.
This separatrix is often referred to as the {\it
strong invariant manifold} of $\fol$. Furthermore there is
also a {\it formal} change of coordinates $(x,y) \mapsto
(\varphi (x,y) ,y)$ which brings $\omega$ to the form
\begin{equation}
\omega (x,y) = x(1 +\lambda y^{p+1}) \, dy \; - y^{p+1} \, dx
\; . \label{fnf}
\end{equation}
The expression in (\ref{fnf}) is said to be the {\it formal
normal form} of $\fol$. In these formal coordinates the
axis $\{ x=0\}$ is invariant by $\fol$ and called the
{\it weak invariant manifold} of $\fol$. Note however that
the weak invariant manifold of $\fol$ does not necessarily
correspond to an actual separatrix of $\fol$ since the change
of coordinates $(x,y) \mapsto (\varphi (x,y) ,y)$ does not
converge in general. Finally it is also known that a saddle-node
$\fol$ possesses at least one and at most two separatrizes
(which are necessarily smooth) depending on whether or not the
weak invariant manifold of $\fol$ is convergent.

A general remark about saddle-nodes is the following one:
denoting by $\pi_2$ the projection $\pi_2 (x,y) = y$,
Dulac's normal form implies
that the fibers of $\pi_2$, namely the vertical lines, are
transverse to the leaves of $\fol$ away from $\{ y=0\}$.
This allows us to define the monodromy of $\fol$ as being
the ``first return map'' to a fixed fiber.

As to semi-complete vector fields whose associated foliation
$\fol$ is a saddle-node, one has:

\begin{teo}
\label{selano}
Suppose that $Y$ is a meromorphic semi-complete vector field
defined around $(0,0) \in \C^2$. Suppose that the foliation
$\fol$ associated to $Y$ is a saddle-node. Then, up to an
inversible factor, $Y$ has one of the following normal forms:
\begin{enumerate}
\item $Y = x(1+ \lambda y) \parlx + y^2 \parly$, $\lambda \in \Z$.

\item $Y = y^{-p} [(x(1+ \lambda y^p) + yR(x,y)) \parlx +
y^{p+1} \parly]$.

\item $Y = y^{1-p} [(x(1+ \lambda y^p) + yR(x,y)) \parlx +
y^{p+1} \parly]$ and the monodromy induced by $\fol$ is trivial
(in particular $\lambda \in \Z$ and the weak invariant manifold
of $\fol$ is convergent).
\end{enumerate}
\end{teo}

\noindent {\it Proof}\,: The proof of this theorem relies heavily
on the methods introduced in Section~4 of~\cite{re2} and Section~4
of~\cite{re4}. For convenience of the reader we summarize the
argument below.

First we set $Y=fZ/g$ where $Z$ is a holomorphic vector field
with an isolated singularity at $(0,0)$. Note that when $f(0,0)
.g(0,0) \neq 0$ (i.e. when $Y$ is holomorphic with an isolated
singularity at $(0,0)$), then $Y$ has the normal form~1. Indeed
this is precisely the content of Theorem~4.1 in Section~4 of
\cite{re4}.

Next we observe that the pole divisor $\{ g=0\}$ is contained
in the strong invariant manifold of $\fol$. To verify this assertion
we first notice that $\{ g=0\}$ must be invariant under $\fol$ as
a consequence of the general fact that there is no one-dimensional
meromorphic semi-complete vector field (cf. Remark~(\ref{obsidiota})).
Thus $\{ g=0\}$ is contained in the union of the separatrizes
of $\fol$. Next we suppose for a contradiction that the weak
invariant manifold of $\fol$ is convergent (i.e. defines a separatrix)
and part of the pole divisor of $Y$. In this case the technique
used in the proof of Proposition~4.2 of \cite{re2} applies
word-by-word to show that the resulting vector field $Y$ is not
semi-complete. This contradiction implies that $\{ g=0 \}$ must
be contained in the strong invariant manifold of $\fol$ as desired.

Combining the information above with Dulac's normal form, we conclude
that $Y$ possesses the form
$$
Y = \frac{f}{y^k} \left[ (x(1+ \lambda y^p) + yR(x,y)) \partx +
y^{p+1} \party \right] \, .
$$
Now we are going to prove that $f(0,0) \neq 0$ i.e. $f$ is an
inversible factor. Hence we assume for a contradiction that
$f(0,0) =0$ but $f$ is not divisible by $y$. Still
according to the terminology of Section~4 of \cite{re2},
we see that the ``asymptotic order of the divided time-form''
induced by $Y$ on $\{ y=0\}$ is at least $2$ since this
form is $dx /(xf(x,0))$ (this is also a consequence of the fact
that the index of the strong invariant manifold of a saddle-node
is {\it zero}, cf. Section~5). However this order cannot be greater
than~$2$ since $Y$ is semi-complete. Furthermore when this order
happens to be~$2$, the local holonomy of the separatrix in question must be
the identity provided that $Y$ is semi-complete. Nonetheless the
local holonomy of the strong invariant manifold of a saddle-node
is never the identity. In fact, using Dulac's normal form,
an elementary calculation shows that this holonomy has
the form $H(z) = z + z^p + \cdots$. We then conclude that
$f(0,0) \neq 0$.

Therefore we have so far
$$
Y = y^{-k} H(x,y) \left[ (x(1+ \lambda y^p) + yR(x,y)) \partx +
y^{p+1} \party \right] \, ,
$$
where $H$ is holomorphic and satisfies $H(0,0) \neq 0$.

Recall that $\pi_2$ denotes the projection $\pi_2 (x,y) =y$
whose fibers are transverse to the leaves of $\fol$ away
from $\{ y=0\}$. Let $L$ be a regular leaf of $\fol$ and
consider an embedded curve $c: [0,1] \rightarrow L$. If $dT_L$
stands for the time-form induced on $L$ by $Y$, we clearly have
$$
\int_c \, dT_l \; = \; \int_{\pi_2 (c)} \, h(y) y^{p-k -1} dy \, ,
$$
where $h (y) = H(0,y)$ so that
$h (0) \neq 0$. Since the integral on the left hand
side is never {\it zero}, it follows that $p-k -1 =0$ or~$1$. The case
$k=p-1$ does not require further comments. On the other hand,
if $k=p-2$, then the integral of $dT_L$ over $c$ is {\it zero}
provided that $\pi_2 (c)$ is a loop around the origin $0 \in
\{ x=0\}$. This implies that $c$ must be closed itself.
In other words the monodromy of $\fol$ with respect to the fibration
induced by $\pi_2$ is trivial. Conversely it is easy to check
that the normal forms 1, 2 and 3 in the statement are, in fact,
semi-complete. This finishes the proof
of the theorem.\qed

Before closing the section, it is interesting to translate the
condition in the item~3 of Theorem~(\ref{selano}) in terms
of the classifying space of Martinet-Ramis \cite{mara}. Note
however that this translation will not be needed for
the rest of the paper.

Fix $p \in \N^{\ast}$ and consider the foliation $\fol_{p, \lambda}$
whose leaves are ``graphs'' (over
the $y$-axis) of the form
$$
x = {\rm const}\, .\,  y^{\lambda} \exp (-1/px^p) \, .
$$
Given $\lambda \in \C$, the moduli space of saddle-nodes $\fol$
having $p, \lambda$ fixed is
obtained from the foliation above through the
following data:
\begin{itemize}

\item $p$ translations $z \mapsto z+ c_i$, $z, c_i \in \C$
denoted by $g_1^+ , \ldots , g_p^+$.

\item $p$ local diffeomorphisms $z \mapsto z + \cdots$, $z \in \C$,
tangent to the identity denoted by $g_1^- , \ldots , g_p^-$.
\end{itemize}

These diffeomorphisms induce a permutation of (part of)
the leaves of $\fol_{p, \lambda}$. More precisely the total
permutation (after one tour around $0 \in \C$) is given by the
composition
$$
g_p^- \circ g_p^+ \circ \ldots \circ g_1^- \circ g_1^+ \, .
$$
However recall that our saddle-node has trivial monodromy.
One easily checks that this cannot happen in the presence of
the ramification $y^{\lambda}$. Thus we conclude that $\lambda$
belongs to $\Z$ and, in particular, the model
$\fol_{p, \lambda}$ introduced above has itself trivial
monodromy. Hence the monodromy of $\fol$ itself is nothing but
$g_p^- \circ g_p^+ \circ \ldots \circ g_1^- \circ g_1^+$.
In other words the condition is
$$
g_p^- \circ g_p^+ \circ \ldots \circ g_1^- \circ g_1^+ = Id \, .
$$
In particular note that, if $p =1$, the above equation
implies that $g_1^- = g_1^+ = Id$. This explains why in item~1
of Theorem~(\ref{selano}) the saddle-node in question is analytically
conjugate to its formal normal form.

\section{Polynomial vector fields and first integrals}

\hspace{0.4cm} Here we want to specifically consider
complete polynomial vector fields whose associated
foliation $\fol$ has a singularity in the line at
infinity which admits
infinitely many separatrizes. Recall
that such a singularity is said to be
dicritical.
The main result of the section is Proposition~(\ref{jouano})
below.

\begin{prop}
\label{jouano}
Let $X$ be a complete polynomial vector field on $\C^2$ and
let $\fol_X$ denote the foliation induced by $X$ on $\pp$.
Assume that $\fol_X$ possesses a dicritical singularity $p$, belonging
to the line at infinity $\Delta$. Then $\fol_X$ has a meromorphic first
integral on $\pp$. Furthermore, modulo a normalization, the closure of the
regular leaves of $\fol_X$ is isomorphic to $\C \Pp (1)$.
\end{prop}

We begin with a weakened version of this proposition
which is the following lemma.

\begin{lema}
\label{prejoa}
Let $X ,\; \fol_X$ be as in the statement of
Proposition~(\ref{jouano}) and denote by $p \in \Delta$
a dicritical singularity of $X$. Then $X$ possesses a meromorphic
first integral on $\C^2$. Moreover, if this first integral
is not algebraic, then the set of leaves of $\fol_X$ that pass through
$p$ only once contains an open set.
\end{lema}

\noindent {\it Proof}\,: Suppose that $\fol_X$ as above possesses
a singularity $p \in \Delta$ with infinitely many separatrizes.
We consider coordinates $(u,v)$ around $p$ such that $\{ u=0 \}
\subset \Delta$. Given a small neighborhood $V$ of $p$, we consider
the restriction $X_{\mid_V}$ of $X$ to $V$. Clearly $X_{\mid_V}$ defines
a meromorphic semi-complete vector field on $V$.

Obviously only a finite number of separatrizes of $\fol_X$
going through $p$ may be contained in the divisor of poles
or zeros of $X$. Thus there are (infinitely many) separatrizes
$\cals$ of $\fol_X$ at $p$ which are (local) regular orbits
of $X$. We fix one of these separatrizes $\cals$. Recall that
$\cals$ has a Puiseux parametrization $\textsf{A}(t) = (a(t), b(t))$,
$\textsf{A}(0) = (0,0)$,
defined on a neighborhood $W$ of $0 \in \C$. Furthermore
$\textsf{A}$ is injective on $W$ and a diffeomorphism from
$W \setminus \{ 0 \}$ onto $\cals \setminus \{ (0,0) \}$.
Since $\cals$ is invariant under $X$, the restriction
to $\cals \setminus \{ (0,0) \}$ of $X$ can be pulled-back
by $\textsf{A}$ to give a meromorphic vector field $Z$ on $W$, i.e.
$Z(t) = \textsf{A}^{\ast} \left( X\! \! \mid_{\cals} \right)$
where $t \in W \setminus \{ (0,0) \}$ for a sufficiently
small neighborhood $W$ of $0 \in \C$. We also have that
$Z$ is semi-complete on $W$ since $X_{\mid_V}$ is semi-complete
on $V$ and $\textsf{A}$ is injective on $W$. It follows from
Remark~(\ref{obsidiota}) that $Z$ admits a holomorphic
extension to $0 \in \C$ which is still denoted by $Z$.
Moreover, letting $Z(t) = h(t) \partial /\partial t$,
we cannot have $h(0) = h'(0) = h'' (0) =0$. However we
must have at least $h(0) =0$. Otherwise the semi-global flow
of $Z$ would reach the origin $0 \in \C$ in finite time.
Since $\textsf{A} (0)$ lies in the line at infinity $\Delta$, it would follow
that points in the orbit of $X$ containing $\cals$ reach $\Delta$
in finite time. This is impossible since $X$ is complete on $\C^2$.

Therefore we have only two cases left, namely $h(0) =0$ but
$h'(0) \neq 0$ and $h(0) = h'(0) =0$ but $h''(0) \neq 0$.
Let us discuss them separately. First suppose that
$h(0) = h'(0) =0$ but
$h'' (0) \neq 0$. Modulo a normalization we can suppose that
$\cals$ is smooth at $p$. Again we denote by $L$ the global
orbit of $X$ containing $\cals$. By virtue of the preceding
$L$ is a Riemann surface endowed with a complete holomorphic vector
field $X_{\mid_L}$ which has a singularity of order~$2$ (where
$X_{\mid_L}$ stands for the restriction of $X$ to $L$).
It immediately results that $L$ has to be compactified into
$\C \Pp (1)$. In other words the closure of $L$ is a rational
curve and, in particular, an algebraic
invariant curve of $\fol_X$. Since there are infinitely many
such curves, Jouanolou's theorem \cite{joa} ensures that
$\fol_X$ has a meromorphic first integral (alternatively
we can also apply Borel-Nishino's theorem, cf. \cite{la}). Furthermore the
level curves of this first integral are necessarily rational
curves (up to normalization) as it follows from the discussion
above.

Suppose now that $h(0) =0$
but $h'(0) \neq 0$. In this case the vector field $Z$
has a non-vanishing residue at $0 \in \C$. We then conclude
that $\cals$ possesses a {\it period} i.e. there exists
a loop $c:[0,1] \rightarrow \cals$ ($c(0) = c(1)$) on which
the integral of the corresponding time-form is different from
{\it zero}. If $L$ is the global orbit of $X$ containing $\cals$
the preceding implies that $L$ is isomorphic to $\C^{\ast}$.

On the other hand the characterization of singularities with
infinitely many separatrizes obtained through Seidenberg's theorem
\cite{sei} (cf. Section~7 for further details) ensures that the
set of orbits $L$ as above has positive logarithmic capacity.
In fact it contains an open set. Thus Suzuki's results in
\cite{suzu1}, \cite{suzu2} apply to provide the existence
of a non-constant meromorphic first integral for $X$.
If this
first integral is not algebraic, then a ``generic'' orbit
passes through $p$ once but not twice. Otherwise the leaf would
contain two singularities of $X$ and therefore it would be a rational
curve (up to normalization). In this case the mentioned Jouanolou's
theorem would provide an algebraic first integral for
$\fol_X$. The proof of the proposition is over.\qed

The preceding lemma suggests a natural strategy to establish
Proposition~(\ref{jouano}). Namely we assume for a contradiction
that $\fol_X$ does not have an algebraic first integral. Then we
have to show that a generic
leaf passing through a dicritical singularity $p$ must return
and cross $\Delta$ once again (maybe through another dicritical
singularity). The resulting contradiction will
then complete our proof.

Using Seidenberg's theorem, we reduce the singularities of $\fol_X$
in $\Delta$ so that they all will have at least one eigenvalue
different from {\it zero}. In particular we obtain a normal crossing
divisor $\cale$ whose irreducible components are rational curves, one of them
being the proper transform of $\Delta$ (which will still be denoted
by $\Delta$). The other components were introduced by the punctual
blow-ups performed and are denoted by $D_i$, $i =1, \ldots ,s$.
The fact that $p$ is a dicritical singularity implies that one
of the following assertions necessarily holds:
\begin{enumerate}

\item There is a component $D_{i_0}$ of $\cale$ which is not
invariant by $\wifol_X$ (where $\wifol_X$ stands for the proper
transform of $\fol_X$).

\item There is a singularity $p_0$ of $\wifol$ in $\cale$
which is dicritical and has two eigenvalues different from
{\it zero}.
\end{enumerate}

We fix a local cross section $\Sigma$ through a point $q$
of $\Delta$ which is regular for $\wifol_X$. Note that a regular leaf $L$ of
$\wifol_X$ necessarily meets $\Sigma$ infinitely many times
unless $L$ is algebraic. Indeed first we
observe that
$L$ is properly embedded in the affine part $\C^2$ since $\fol_X$ possesses
a meromorphic first integral on $\C^2$. Thus all the accumulation points
of $L$ are contained in $\Delta$. Obviously if $L$ accumulates a regular
point of $\Delta$ then $L$ intersects $\Sigma$ infinitely many times
as required. On the other hand, if $L$ accumulates only points of
$\Delta$ which are singularities of $\fol_X$, then Remmert-Stein theorem
shows that the closure of $L$ is algebraic. Summarizing, using Jouanolou's
or Borel-Nishino's theorem, we can suppose without loss of generality
that all the leaves of $\wifol_X$ intersects $\Sigma$ an infinite number
of times (and in fact these intersection points approximate
the point $q = \Sigma \cap \Delta$).

To prove that $\fol_X$ has infinitely many leaves cutting
the exceptional divisor $\cale$ more than once, we fix a neighborhood
$\mathcal U$ of $\cale$. Proposition~(\ref{jouano}) is now
a consequence of the next proposition.

\begin{prop}
\label{dutrans}
Under the above assumptions, there is an open neighborhood $V
\subset \Sigma$ of $q$ in $\Sigma$ and an open subset
$W \subset V$ of $V$ with the following property: any leaf $L$
passing through a point of $W$ intersects the exceptional
divisor $\cale$ before leaving the neighborhood $\mathcal U$.
\end{prop}

In fact, since the set of leaves meeting $\cale$ contains an
open set and all of them (with possible exception of a finite number)
cross $\Sigma$ and accumulates $\Delta$,
Proposition~(\ref{dutrans}) clearly shows the existence of infinitely
many leaves (orbits of $X$) intersecting $\cale$ more than
one time thus providing the desired contradiction.

In order to prove Proposition~(\ref{dutrans}) we keep the preceding
setting and notations. We are naturally led to discuss the
behavior of the leaves of $\fol_X$ on a neighborhood of the
point $p_{ij}$ of intersection of the irreducible components
$D_i , D_j$ belonging to $\cale$.

Now let us fix coordinates $(x,y)$ around $p_{ij}$ ($p_{ij}
\simeq (0,0)$) such that $\{ y= 0\} \subseteq D_i$ and
$\{ x=0 \} \subseteq D_j$. Without loss of generality we can
suppose that the domain of definition of the $(x,y)$-coordinates
contains the bidisc of radius~$2$. Next we fix a segment of vertical line
$\Sigma_x$ (resp. horizontal line $\Sigma_y$) passing through
the point $(1,0)$ (resp. $(0,1)$). We assume that $D_i$ is
invariant under $\wifol_X$ but $D_j$ may or may not be invariant
under $\wifol_X$. Let us also make the following assumptions:
\begin{description}

\item[{\bf A)}] $p_{ij} \simeq (0,0)$ is not a dicritical
singularity of $\wifol$.

\item[{\bf B)}] $\wifol_X$ has at least one eigenvalue different
from {\it zero} at $p_{ij} \simeq (0,0)$.

\item[{\bf C)}] The vector field $X$ whose associated foliation
is $\wifol$ is meromorphic semi-complete in the domain of the
coordinates~$(x,y)$.
\end{description}

We are going to discuss a variant of the so-called {\it Dulac's transform},
namely if leaves intersecting $\Sigma_x$ necessarily cut $\Sigma_y$.
Precisely we fix a neighborhood ${\bf U}$ of $\{ x=0 \} \cup
\{ y=0 \}$, we then have:

\begin{lema}
\label{passara}
Under the preceding assumptions, there is an open
neighborhood $V_x \subset \Sigma_x$ of $(1,0)$ in
$\Sigma_x$ and an open set $W_x \subset V_x$
with the following property: any leaf $L$ of $\fol$
passing through a point of $W_x$ meets $\Sigma_y$
before leaving $U$. In particular, if $D_j$ is not
invariant by $\fol$, then the leaves of $\fol$ passing
through points of $W_x$ cross the axis~$\{ x=0 \}$
before leaving ${\bf U}$.
In addition, by choosing $V_x$ very small, the ratio
between the area of $W_x$ and the area of $V_x$ becomes
arbitrarily close to~$1$.
\end{lema}

Before proving this lemma, let us deduce the proof
of Proposition~(\ref{dutrans}).

\vspace{0.1cm}

\noindent {\it Proof of Proposition~(\ref{dutrans})}\,: Recall
that $\wifol_X$, the proper transform of $\fol_X$, has only
simple singularities in $\cale$. In particular, if ${\bf p}
\in \cale$ is a dicritical singularity of $\wifol_X$, then $\wifol_X$
has~$2$ eigenvalues different from {\it zero} at ${\bf p}$.
Hence $\wifol_X$ is linearizable around ${\bf p}$ and, as a
consequence, there is a small neighborhood $U_{\bf p}$ of
${\bf p}$ such that any regular leaf $L$ of $\wifol_X$ entering
$U_{\bf p}$ must cross $\cale$ before leaving
$U_{\bf p}$. This applies in particular if ${\bf p}$ coincides
with the intersection of two irreducible components $D_i,
D_j$ of $\cale$.

On the other hand $\Delta$ is invariant by $\wifol_X$ and
all but a finite number of leaves of $\wifol_X$ intersect $\Sigma$
in arbitrarily small neighborhoods of $q= \Sigma \cap \Delta$.
In particular if $\wifol_X$ has a dicritical singularity on
$\Delta$, then the statement follows from the argument above.

Suppose now that $\wifol_X$ does not have a dicritical singularity
on $\Delta$. Let $D_1$ be another irreducible component of
$\cale$ which intersects $\Delta$ at $p_{01}$. Note that
$p_{01}$ is not a dicritical singularity of $\wifol_X$ by the
preceding discussion. If $D_1$ is not invariant by $\wifol_X$,
then Lemma~(\ref{passara}) allows us to find infinitely many
leaves of $\wifol_X$ intersecting $\cale$. Thus the proposition
would follow. On the other hand, if $D_1$ is invariant by
$\wifol_X$, then Lemma~(\ref{passara}) still allows us to find
a local transverse section $\Sigma_1$ through a regular point
$q_1$ of $D_1$ with the desired property namely: excepted for a 
set of leaves whose volume can be made arbitrarily small
(modulo choosing $V_x$ sufficiently small), all leaves of $\wifol_X$
meet $\Sigma_1$ in arbitrarily small neighborhoods of
$q_1 = \Sigma_1 \cap D_1$. We then continue the procedure replacing
$\Delta$ by $D_1$. Since we eventually will find an irreducible
component of $\cale$ which is not invariant by $\wifol_X$ or contains
a dicritical singularity, the proposition is proved. This also
concludes the proof of Proposition~(\ref{jouano}).\qed

The rest of the section is devoted to the proof of
Lemma~(\ref{passara}). Let us begin with the easier case
in which $D_j$ is not invariant by $\wifol_X$.

\begin{lema}
\label{10.fim1}
The vector field $\widxx$ vanishes with order~$1$ on $D_j$.
It also has poles of order~$1$ on $D_i$ and the origin
$(0,0)$ is a LJ-singularity of $\wifol_X$.
\end{lema}

\noindent {\it Proof}\,: By assumption $D_j$ is contained in
$\cale$ and is not invariant by $\wifol_X$. In particular
$\widxx$ cannot have poles on $D_j$ since there is no
strictly meromorphic semi-complete vector field in dimension~$1$.
Neither can $\widxx$ be regular on $D_j$ otherwise certain points
of $\C^2$ would reach the infinity in finite time, thus contradicting
the fact that $\widxx$ is complete. Finally the order of $\widxx$
on $D_j$ cannot be greater than~$2$, otherwise $\widxx$
would not be semi-complete. Besides, if this order is~$2$, then
infinitely many orbits of $X$ will be compactified into rational
curves. This would imply that $\fol_X$ has an algebraic first
integral which is impossible. This shows that $\widxx$ vanishes
with order~$1$ on $D_j$.

Because $\widxx$ vanishes on $D_j$ and $D_j$ is not
invariant by the associated foliation $\wifol_X$, it follows
from Section~4 that either $(0,0)$ is a LJ-singularity of
$\wifol_X$ or $\wifol_X$ is linearizable with eigenvalues~$1,-1$
at $(0,0)$. In the latter case the conclusions of
Lemma~(\ref{passara}) are obvious. Thus we can suppose that
$(0,0)$ is a LJ-singularity. It follows from Lemma~(\ref{le3.3})
that $\widxx$ must have a pole divisor of order~$1$ on $D_i$.\qed

\vspace{0.1cm}

\noindent {\it Proof of Lemma~(\ref{passara}) when $D_j$
is not invariant by $\wifol_X$}\,: Modulo blowing-up $(0,0)$,
the problem is immediately reduced to the discussion of the Dulac's
transform between the strong and the weak invariant
manifolds of the saddle-node determined by
$$
x(1+y) \, dy \; \,  + \; \,  y^2 \, dx \; .
$$
Thinking of this foliation as a differential equation,
we obtain the solution
$$
x(T)  =  \frac{x_0 e^T}{1-y_0 T} \; \; \, {\rm and} \; \; \,
y(T)  =  \frac{y_0}{1-y_0 T} \, .
$$
Let us fix $x_0 =1$. Given $y_0$, we search for $T_0$ so that
$y(T_0) =1$. Furthermore we also require that the norm of $x(T)$
stays ``small'' during
the procedure. This is clearly possible if $y_0$ is real negative
(sufficiently close to {\it zero}). Actually we set
$T_0 =(1-y_0) /y_0 \in \R_-$ so that $T$ can be chosen real negative during
the procedure. Thus the norm of $x(T)$ will remain controlled by
that of $y_0$.

On the other hand, let $y_0$
be in the transverse section $\Sigma_x$ and suppose that $y_0$
is not real positive. Then the orbit of $y_0$ under the local
holonomy of the strong invariant manifold converges to {\it zero}
and is asymptotic to $\R_-$. Indeed this local holonomy is represented
by a local diffeomorphism of the form $z \mapsto z + z^2 +
{\rm h.o.t.}$ and the local topological dynamics of these
diffeomorphisms is simple and well-understood (known as a
``flower'', in the present case this dynamics is also called
the parabolic bifurcation). Hence for a sufficiently large iterate
of $y_0$ (without leaving $V_x$), the above ``Dulac's transform'' is well-defined.
Thus the statement of Lemma~(\ref{passara}) is verified
as long as we take $W_x = V_x \setminus \R_-$ in the above
coordinates.\qed

From now on we can suppose that $D_j$ is invariant under $\fol_X$.
We have three cases to check:

\noindent 1) $\fol_X$ has two eigenvalues different from {\it zero}
at $(0,0)$ and is locally linearizable.

\noindent 2) $\fol_X$ has two eigenvalues different from {\it zero}
at $(0,0)$ but is not locally linearizable. In this case
the quotient of the eigenvalues is real negative.

\noindent 3) $\fol_X$ defines a saddle-node at $(0,0)$.

In the Case~$1$ the verification is automatic and left to the reader.
Case~2 follows from \cite{mamo} (note that our convention of signs
is opposite to the convention of \cite{mamo}). So we just need
to consider the case of saddle-nodes. Of course all the
possible saddle-nodes necessarily have a convergent weak invariant
manifold. Without loss of generality we can suppose that
$D_i$ is the strong invariant manifold so that $D_j$ is the
weak invariant manifold (the other possibility is analogous).
All the background material about saddle-nodes
used in what follows can be found in \cite{mara}.

Thanks to Lemma~(\ref{selano}), we can find coordinates $(x,y)$
as above where the vector field $\widxx$ becomes
$$
\widxx = y^{-k} [ (x(1+ \lambda y^p) +yR) \parlx \; + \;
y^{p+1} \parly ] \, .
$$
Since our problem depends only on the foliation associated
to $\widxx$, we drop the factor $y^{-k}$ in the sequel. We also
notice that $\widxx$ is regular on $D_j$. The argument which is
going to be
employed here is a generalization of the one employed to deal with
the saddle-node appearing after blowing-up the LJ-singularity
in the previous case.

Following \cite{mara} we consider open sets $V_i \subset \C$,
$i = 0, \ldots ,2p-1$, defined by $V_i = \{ z \in \C \; ; \;
(2i+1)\pi /2p - \pi /p < \arg z < (2i+1)\pi /2p + \pi /p
\}$. The $V_i$'s, $i=1, \ldots ,2p-1$, define a covering
of $\C^{\ast}$ and, besides, each $V_i$ intersects only
$V_{i-1}$ and $V_{i+1}$ (where $V_{-1} = V_{2p-1}$). We also
let $W_i^+$ (resp. $W_i^-$) be defined by
$$
W_i^+ =  \left\{ \frac{(4i-1) \pi}{2p} < \arg z <
\frac{(4i+1) \pi}{2p} \right\} \, \; \; {\rm and} \, \; \;
W_i^-  =  \left\{ \frac{(4i+1) \pi}{2p} < \arg z <
\frac{(4i+3) \pi}{2p} \right\} \, .
$$
for $i=0 ,\ldots ,2p-1$. We point out that ${\rm Re} \, (y^p)
< 0$ (resp. ${\rm Re} \, (y^p) > 0$) provided that $y \in W_i^-$
(resp. $y \in W_i^+$). In addition we have $W_i^+ = V_{2i}
\cap V_{2i+1}$ and $W_i^- = V_{2i+1} \cap V_{2i+2}$
(unless $p=1$ where $V_0 \cap V_1 = W_0^+ \cup W_0^-$).
Given $\varepsilon >0$, we set
$$
U_{i, V} = \{ (x,y) \in \C^2 \; ; \; \| x \| < \varepsilon
\, , \; \, \| y \| < \varepsilon \; \; {\rm and} \; \;
y \in V_i \} \, .
$$
According to Hukuara-Kimura-Matuda, there is a bounded
holomorphic mapping $\phi_{U_{i, V}} (x,y) = (\varphi_{U_{i,V}} (x,y) , y)$
defined on $U_{i, V}$ which brings the vector field $\widxx$
to the form
\begin{equation}
x(1+ \lambda y^p) \parlx \; + \; y^{p+1} \parly \, .
\label{quasela}
\end{equation}
The vector field in~(\ref{quasela}) can be integrated to
give
\begin{equation}
x (T) = \frac{x_0 e^T}{\sqrt[p]{(1-py_0^pT)^{\lambda}}}
\; \; {\rm and} \; \; y(T) = \frac{y_0}{\sqrt[p]{1-p
y_0^p T}} \, . \label{intee}
\end{equation}

\noindent {\it Proof of Lemma~(\ref{passara})}\,: We keep the
preceding notations. On $U_{i,V}$ we consider a normalizing mapping
$\phi_{U_{i, V}}$ such that $\widxx$ is as in~(\ref{quasela}).
In this coordinate we fix a vertical line
$\Sigma_x$ as before and let $\Sigma_{y,i}$
denote the intersection of the horizontal line through $(0,1)$,
$\Sigma_y$, with the sector $V_i$.
Again we want to know which
leaves of $\wifol$ passing through a point of $\Sigma_x$ will
intersect $\Sigma_y$ as well. Thus starting with $x_0 =1$,
we search for $y (T_0) =1$. Thanks to equations~(\ref{intee}),
it is enough to choose $T_0=(1-y_0^p) /py_0^p$. In particular
$T_0 \in \R_-$ provided that $y_0^p \in \R_-$. The formula
for $x(T)$ in~(\ref{intee}) shows that $x(T)$ remains in the
fixed neighborhood ${\bf U}$ of $\{ x=0 \} \cup \{ y=0 \}$ provided
that we keep $T \in \R_-$ during the procedure and choose
$y_0$ sufficiently small. More generally if the real part
${\rm Re} \, (y_0^p)$ is negative and the quotient between
imaginary and real parts is bounded, then the same argument
applies. In other words, if $y_0$ belongs to a compact subsector
of $W_{[i+1/2]-1}^-$, the set $W_j^-$ contained in $V_i$,
then the Dulac's transform in question is well-defined
modulo choosing $y_0$ uniformly small.

The local holonomy associated to $D_i$ (the strong invariant
manifold) has the form $\textsf{h} (z) = z + z^{p+1} + {\rm h.o.t.}$
The dynamical picture corresponding to this diffeomorphims
is still a ``flower''. However,
in general, it is not true that the orbit of a ``generic'' point
will intersect a fixed $W_{[i+1/2]-1}^-$ since $\textsf{h}$ may have
invariant sectors. However, the above argument can be applied
separately for each $i$. Clearly to each $i$ fixed we have a
different $\Sigma_y \cap V_i$ associated. Nonetheless they are all
equivalent since the weak invariant manifold
of the foliation is convergent. It follows that apart from a
finite number of curves whose union has empty interior,
the leaf through a point of $\Sigma_x$ meets $\Sigma_y$
before leaving the neighborhood ${\bf U}$. As mentioned
the case where $D_i$ is the (convergent) weak invariant manifold
and $D_j$ is the strong invariant manifold is analogous.
The statement follows.\qed

\section{Arrangements of simple singularities}

\hspace{0.4cm}
Now we are going to study the possible arrangements of
simple semi-complete singularities over a rational
curve of self-intersection~$-1$ which are obtained by blowing-up
a semi-complete vector field on a neighborhood of
$(0,0) \in \C^2$.

We shall make a number of assumptions which are always satisfied in
our cases.
We denote by
$\wicc$ the blow-up of $\C^2$ at the origin
and by $\pi : \wicc \rightarrow \C^2$ the corresponding
blow-up map. Given a vector field $Y$ (resp. foliation
$\fol$) defined on a neighborhood $U$ of the origin,
$\pi$ naturally induces a vector field $\wiX$ (resp. foliation $\wifol$)
defined on $\pi^{-1} (U)$. Furthermore $Y$ is semi-complete
on $U$ if and only if $\wiX$ is semi-complete on $\pi^{-1} (U)$.

Given a meromorphic vector field $Y = fZ/g$ with $f,g, Z$
holomorphic, we call the vector field $fZ$ the {\it holomorphic
part}\, of $Y$. In this section it is discussed the nature of a
meromorphic semi-complete vector field $Y$ defined
on a neighborhood of the origin $(0,0) \in \C^2$
which satisfies the following assumptions:
\begin{enumerate}
\item $Y = x^{k_1} y^{k_2} f_1 Z $ where $Z$
is a holomorphic vector field having an isolated singularity
at $(0,0) \in \C^2$, $f_1$ is
holomorphic function and $k_1 , k_2 \in \Z$.

\item The regular orbits of $Y$ contain at most
$1$ singular point. Furthermore the order of $Y$ at
this singular point is one.

\item The regular orbits of $Y$ contains at most~$1$
period (i.e. there is only one homology class
containing loops on which the integral of the time-form
is different from zero).

\item The foliation $\fol$ associated to $Y$ (or to
$Z$) has both eigenvalues equal to {\it zero}
at $(0,0) \in \C^2$.

\item The blow-up $\wifol$ of $\fol$ is such that
every singularity $\tilde{p} \in \pi^{-1} (0)$
of $\wifol$ is simple.

\item $\fol$ is not dicritical at $(0,0)$.

\end{enumerate}

Before continuing let us introduce two basic definitions. Assume
that $\fol$ is a singular holomorphic foliation defined on a
neighborhood of a singular point $p$. Let $\cals$ be a
smooth separatrix of $\fol$ at $p$. We want to define the
{\it order }\, of $\fol$ with respect to $\cals$ at $p$,
${\rm ord}_{\cals} (\fol ,p)$ (also called the multiplicity
of $\fol$ along $\cals$), and
the {\it index}\, of $\cals$ w.r.t. $\fol$ at $p$,
${\rm Ind}_{p} (\fol , \cals)$ (cf. \cite{cama}).
In order to do that, we consider coordinates $(x,y)$ where
$\cals$ is given by $\{ y=0 \}$ and a holomorphic $1$-form
$\omega = F(x,y) \, dy - G (x,y)\, dx$ defining $\fol$ and having
an isolated singularity at $p$. Then we let
\begin{eqnarray}
{\rm ord}_{\cals} (\fol , p) & = & {\rm ord} \, (F(x,0)) \; \; \,
{\rm at} \; \; \, 0 \in \C \label{ordem1} \; \; \, {\rm and} \\
{\rm Ind}_p (\fol ,\cals ) & = & {\rm Res} \,
\frac{\partial}{\partial y} \left( \frac{G}{F} \right) (x,0) \, dx
\; . \label{index1}
\end{eqnarray}
In the above formulas ${\rm ord} \, (F(x,0))$ stands for the order of the
function $x \mapsto F(x,0)$ at $0 \in \C$ and
${\rm Res}$ for the residue of the $1$-form in question.

Let $p_1 ,\ldots ,p_r$ denote the singularities of $\wifol$
belonging to $\pi^{-1} (0)$. Since $\pi^{-1} (0)$ naturally
defines a separatrix for every $p_i$, we can consider both
${\rm ord}_{\pi^{_1} (0)} (\wifol , p_i)$ and ${\rm Ind}_{p_i}
(\wifol ,\pi^{-1} (0))$. Easy calculations and the Residue Theorem
then provides (cf. \cite{mamo} , \cite{cama}):
\begin{eqnarray}
& & {\rm ord}_{(0,0)} (\fol) + 1  =  \sum_{i=1}^r
{\rm ord}_{\pi^{-1} (0)} (\wifol , p_i) \; , \label{ordem2} \\
& & \sum_{i=1}^r {\rm Ind}_{p_i} (\wifol, \pi^{-1} (0)) =
-1 \; . \label{index2}
\end{eqnarray}

On the other hand the order of $\pi^{-1} (0)$ as a divisor
of zeros or poles of $\wiX$ is
\begin{equation}
{\rm ord}_{\pi^{-1} (0)} \wiX = {\rm ord}_{(0,0)} (f)
+ {\rm ord}_{(0,0)} (\fol) -{\rm ord}_{(0,0)} (g) -1  \, . \label{equa2}
\end{equation}
In particular if this order is {\it zero}, then $\wiX$
is regular on $\pi^{-1} (0)$.

To abridge notations, the local singular foliation induced
by the linear vector field
$$
(x+y) \parlx + y \parly
$$
will be called a LJ-singularity (where LJ stands for linear and
in the Jordan form).

To simplify the statement of the main result of this section,
namely Proposition~(\ref{prop4.2}), we first introduce~$3$ types,
or models, of vector fields. Let us keep the preceding notation.
\begin{description}

\item[{\bf Model} $Z_{1,11}$:] Let $\fol_{1,11}$ be the foliation associated
to $Z_{1,11}$ and $\wifol_{1,11}$ its blow-up. Then $\wifol_{1,11}$
contains~$3$ singularities $p_1 ,p_2 ,p_3$ on $\pi^{-1} (0)$ whose
eigenvalues are
respectively $1,1$, $-1,1$ and $-1,1$. The singularity $p_1$ is a
LJ-singularity and the blow-up $\widetilde{Z}_{1,11}$ has a pole
of order~$1$ on $\pi^{-1} (0)$. The separatrix of $p_2$ (resp. $p_3$)
transverse $\pi^{-1} (0)$ is a pole divisor of
$\widetilde{Z}_{1,11}$ of order~$1$ as well. Finally $\widetilde{Z}_{1,11}$
has a curve of zeros passing through $p_1$ which is not invariant
under $\wifol_{1,11}$.

\item[{\bf Model} $Z_{0, 12}$:] With similar notations, $\wifol
_{0,12}$ has~$3$ singularities $p_1, p_2 ,p_3$ on $\pi^{-1} (0)$
of eigenvalues equal to $1,0$, $-1,2$ and $-1, 2$. The singularity
$p_1$ is a saddle-node with strong invariant manifold contained in
$\pi^{-1} (0)$. The separatrix of $p_2$ (resp. $p_3$)
transverse to $\pi^{-1} (0)$ is a pole of order~$d\neq 0$ of
$\widetilde{Z}_{0,12}$. The $\pi^{-1} (0)$ is a pole
of order~$2d-1$ of $\widetilde{Z}_{0,12}$. There is no other component
of the divisor of zeros or poles of $\widetilde{Z}_{0,12}$. 

\item[{\bf Model} $Z_{1, 00}$:] $\wifol_{1,00}$ still has~$3$ singularities
$p_1, p_2 ,p_3$ whose eigenvalues are $-1,1$, $1,0$ and~$1,0$. The
singularities $p_2, p_3$ are saddle-nodes with strong invariant
manifolds contained in $\pi^{-1} (0)$. The separatrix
of $p_1$ transverse to $\pi^{-1} (0)$ is a pole of
$\widetilde{Z}_{1,00}$ of order~$d\neq 0$. The exceptional divisor $\pi^{-1} (0)$
is a pole of order~$d-1$ and there is no other component of the divisor
of zeros or poles of $\widetilde{Z}_{1,00}$,

\end{description}

Note in particular that Formula~(\ref{ordem2}) implies that
$\fol_{1,11}$ (resp. $\fol_{0,12}$, $\fol_{1,00}$) has a singularity
of order~$2$ at the origin.

\begin{prop}
\label{prop4.2}
Let $Y$, $\wiX$ be as above. Assume that the order of
$\wiX$ on $\pi^{-1} (0)$ is different from {\it zero}.
Then the structure of the singularities
of $\wifol$ on $\pi^{-1} (0)$ is equal to that of one of the models
$Z_{1,11}$, $Z_{0,12}$ or $Z_{1,00}$.
\end{prop}

\begin{lema}
\label{lema4.3}
Denote by $\lambda_1^i , \lambda_2^i$ the eigenvalues of
$\wifol$ at $p_i$ ($i=1 ,\ldots ,r$). Then one of the
following possibilities holds:
\begin{description}

\item[($\imath$)] $\lambda_1^i /\lambda_2^i = -n/m$ where
$n,m$ belong to $\N^{\ast}$.

\item [($\imath \imath$)] $p_i$ is a saddle-node (i.e.
the eigenvalues are $1$ and {\it zero}) whose strong
invariant manifold coincides with $\pi^{-1} (0)$.

\item[($\imath \imath \imath $)] $p_i$ is a LJ-singularity.
\end{description}
Furthermore there may exist at most one LJ-singularity and,
when such singularity does exist, all the remaining
singularities are as in ($\imath$).
\end{lema}

\noindent {\it Proof}\,: First let us suppose that
one of the eigenvalues $\lambda_1^i ,\lambda_2^i$
vanishes. In this case $\wifol$ defines a saddle-node
at $p_i$. Moreover $p_i$ belongs to the divisor of zeros or poles
of $\wiX$, so that Theorem~(\ref{selano})
shows that the strong invariant manifold of $\wifol$
at $p_i$ coincides with $\pi^{-1} (0)$ as required.

On the other hand, if both $\lambda_1^i ,\lambda_2^i$ are
different from {\it zero}, then they satisfy condition
($\imath$) as a consequence of Lemma~(\ref{le3.1}).
Finally it remains only to consider the case where
$p_1$ is a LJ-singularity. Thus in local coordinates
$(x,t)$, $\{ x=0 \} \subset \pi^{-1} (0)$, around $p_1$,
$\wiX$ is given by
$$
x^{-1} h [ (t+x) \partial /\partial t + x \parlx] \, .
$$
Hence the regular orbits of $\wiX$ contain a {\it zero}
of $\wiX$ corresponding to their intersection with $h=0$.
Because of condition~2, this implies that only one of the
$p_i$'s can be a LJ-singularity. Furthermore, by the
same reason, the holonomy of $\pi^{-1} (0) \setminus
\{ p_1 ,\ldots ,p_r \}$ has to be trivial. In particular
none of the remaining singularities can be a saddle-node.\qed

Combining the information contained in the preceding
lemma with Formula~(\ref{ordem2}) we obtain:

\begin{coro}
\label{lema4.1}
The order of $\fol$ at $(0,0) \in \C$ equals $r-1$, i.e.
${\rm ord}_{(0,0)} (\fol) = r-1$.\qed
\end{coro}

The case where $p_1$ is a LJ-singularity is indeed easy
to analyse. After the preceding lemma and the fact that
the holonomy of $\pi^{-1} (0) \setminus
\{ p_1 ,\ldots ,p_r \}$ is trivial, we conclude that
all the remaining singularities $p_2 ,\ldots ,p_r$ have
eigenvalues $1$ and $-1$. Now using formulas~(\ref{ordem2})
and~(\ref{index2}) we conclude that $\wiX$ is as in the
model~$Z_{1,11}$.

Hereafter we suppose without loss of generality that none
of the $p_i$'s is a LJ-singularity.
For $s \leq r$, we denote by $p_1 ,\ldots ,p_r$ the
singularities of $\wifol$ where $\wifol$ has two
non-vanishing eigenvalues (whose quotient has the
form $-n/m$, $m,n \in \N$). The remaining $p_{s+1} ,
\ldots ,p_r$ singularities are therefore saddle-nodes.
Recall that
the strong invariant manifolds of these saddle-nodes
coincide with $\pi^{-1} (0)$ thanks to Theorem~(\ref{selano}).
Next we have:

\begin{lema}
\label{lema4.2}
At least one of the $p_i$'s is a saddle-node (i.e.
$s$ is strictly less than $r$).
\end{lema}

\noindent {\it Proof}\,: The proof relies on Section~4
of \cite{re3}. Suppose for a contradiction that none of
the $p_i$'s is a saddle-node. Given that there is no
$LJ$-singularity, it follows that the quotient
$\lambda_1^i /\lambda_2^i$ is negative rational for every
$i=1, \ldots ,r$. Hence the hypotheses of Proposition~4.2
of \cite{re3} are verified. It results that $X$ has
one of the normal forms indicated in that proposition.
As is easily seen, all those vector fields have orbits
with~$2$ distinct periods which is impossible in our
case. The lemma is proved.\qed

To complete the proof of Proposition~(\ref{prop4.2}) we
proceed as follows. For $i \in \{ 1, \ldots ,s\}$, we
consider local coordinates $(x_i ,t_i)$, $\{ x_i =0 \}
\subset \pi^{-1} (0)$, around $p_i$. In these coordinates
$\wiX$ has the form
\begin{equation}
\wiX = x_i^{({\rm ord}_{\pi^{-1} (0)} (\wiX))} t_i^{d_i}
h_i [ m_i x_i (1 + {\rm h.o.t.}) \partial /\partial x_i - n_i t_i
(1 + {\rm h.o.t}) \partial /\partial t_i] \label{forca1}
\end{equation}
where $d_i \in \Z$, $m_i ,n_i \in \N$ and $h_i$ is holomorphic
but not divisible by either $x_i, t_i$. Similarly, around
the saddle-nodes singularities $p_{s+1} ,\ldots ,p_r$, we have
\begin{equation}
\wiX = x_i^{({\rm ord}_{\pi^{-1} (0)} (\wiX))} h_i[x_i^{p_i+1}
\partial /\partial x_i + t_i(1 + {\rm h.o.t}) \partial
/\partial t_i] \, . \label{forca2}
\end{equation}

We claim that $h_i (0,0) \neq 0$. This is clear in equation~(\ref{forca2})
thanks to Theorem~(\ref{selano}). As to equation~(\ref{forca1}), let
us suppose for
a contradiction that $h_i (0,0) =0$. Hence the regular leaves of
$\wiX$ have a zero corresponding to the intersection of these
leaves with $\{ h_i =0 \}$. Given condition~2, it results
that only one of the $h_i$'s may verify $h_i (0,0) =0$.
Without loss of generality we suppose that $h_1 (0,0) =0$.
From Lemma~(\ref{le3.1}) and Remark~(\ref{ipc}), it follows
that $(x_i ,t_i)$ can be chosen so as to have $\wiX =
(xy)^a (x-y) (x_i \partial /\partial x_i - t_i \partial /\partial
t_i)$. Formula~(\ref{index2}) then shows that all the
remaining singularities have to be saddle-nodes since the
sum of the indices is~$-1$. Nonetheless, again condition~2, implies
that the holonomy of $\pi^{-1} (0) \setminus \{ p_1 ,\ldots 
, p_r \}$ is trivial. Thus no saddle-node can appear on $\pi^{-1}
(0)$. In other words $r$ must be equal to~$1$ which is impossible.

\vspace{0.1cm}

\noindent {\it Proof of Proposition(\ref{prop4.2})}\,: Considering
the normal forms~(\ref{forca1}) and~(\ref{forca2}), we can suppose
that $h_i (0,0) \neq 0$. Set $\epsilon_i = ({\rm ord_{\pi^{-1} (0)}}
(\wiX)) m_i - n_i d_i$ so that $\epsilon_i \in \{ -1, 0 ,1\}$
(cf. Lemma~(\ref{le3.1})). Alternatively we let $d_i=
({\rm ord_{\pi^{-1} (0)}} (\wiX)) m_i/n_i - \epsilon_i /n_i$.

On the other hand, Formula~(\ref{index2}), in the present
context, becomes
$$
\sum_{i=1}^r m_i /n_i = 1 \, ,
$$
where $m_i=0$ if and only if $p_i$ is a saddle-node
and $n_i \neq 0$. Since all $m_i ,n_i$ are non-negative,
only one of the $n_i$'s can be equal to~$1$ provided that
$m_i \neq 0$. In this case, we must have $m_i=1$ as well
and the remaining singularities are saddle-nodes. We claim that
this implies that $h_i (0,0)$ in~(\ref{forca1}) is always
different from {\it zero}. Indeed if, say $h_1 (0,0) =0$,
then $m+1 =n_1 =1$ and the remaining singularities are saddle-nodes.
The fact that the holonomy associated to the strong
invariant manifold of a saddle-node is has order infinity,
implies that this case cannot be produced. The resulting
contradiction establishes the claim.

Now the fact that $h_i (0,0) \neq 0$ show that
$\sum_{i=1}^r d_i = {\rm ord}_{(0,0)} (f) -
{\rm ord}_{(0,0)} (g)$. Therefore
\begin{eqnarray*}
{\rm ord}_{(0,0)} (f) -
{\rm ord}_{(0,0)} & = & \sum_{i=1}^r d_i =
({\rm ord_{\pi^{-1} (0)}} (\wiX)) (\sum_{i=1}^r m_i /n_i = 1)
- \sum_{i=1}^r \epsilon_i /n_i \\
 & = & - \sum_{i=1}^r \epsilon_i /n_i +
{\rm ord}_{(0,0)} (f) + r-1 -{\rm ord}_{(0,0)} (g) -1 \, .
\end{eqnarray*} 
In other words, one has
\begin{equation}
\sum_{i=1}^s (1- \epsilon_i /n_i) = 2 + s -r <2\, .\label{2sr=2}
\end{equation}
As mentioned, only one of the $n_i$'s may be equal to~$1$.
In this case the remaining singularities are
saddle-nodes and we obtain the model~$Z_{1,00}$.

Next assume that all the $n_i$'s are strictly greater than~$1$.
In particular $1-\epsilon_i /n_i \geq 1/2$. The only new possibility
is to have $n_1=n_2=2$ and $r-s=1$. Thus we obtain the model~$Z_{0,12}$
completing the proof of our proposition.\qed

\begin{obs}
\label{locallinear}
{\rm To complement the description of the Models $Z_{1,11}$,
$Z_{0,12}$ and $Z_{1,00}$, we want to point out that excepted for
the saddle-nodes, all
the singularities appearing in the exceptional divisor after
blowing-up are linearizable. Indeed this results from the
finiteness of the local holonomies associated to their
separatrizes. To check that these holonomies are finite
we just have to use an argument analogous to the one employed
in Remark~(\ref{ipc}).

As a consequence of the above fact, we conclude that the two
saddle-nodes appearing as singularities of $\wifol_{1,00}$
are identical. In particular either both have convergent weak
invariant manifold or both have divergent weak invariant manifold.
}
\end{obs}

\section{The combinatorics of the reduction of singularities}

\hspace{0.4cm}
In this last section we are going to prove our main results.
Since we are going to work in local coordinates, we can consider
a meromorphic semi-complete vector field $Y$ defined around
$(0,0) \in \C^2$. As usual let $\fol$ be the foliation associated
to $Y$. In view of Seidenberg's theorem \cite{sei}, there exists
a sequence of punctual blow-ups $\pi_j$ together with
singular foliations $\wifol^j$,
\begin{equation}
\fol = \wifol^0 \stackrel{\pi_1}{\longleftarrow} \wifol^1
\stackrel{\pi_2}{\longleftarrow} \cdots
\stackrel{\pi_r}{\longleftarrow} \wifol^r \, , \label{resotree}
\end{equation}
where $\wifol^j$ is the blow-up of $\wifol^{j-1}$, such that
all singularities of $\wifol^r$ are simple. Furthermore each
$\pi_j$ is centered at a singular point where $\wifol^{j-1}$
has vanishing eigenvalues. The sequence $(\wifol^j , \pi_j)$ is
said to be the {\it resolution tree} of $\fol$. Fixed
$j \in \{ 1, \ldots ,r \}$, we denote by $\cale^j$
the total exceptional divisor $(\pi_1 \circ \cdots \circ
\pi_j)^{-1} (0,0)$ and by $D^j$ the irreducible component
of $\cale^j$ introduced by $\pi_j$. Note that $D^j$ is
a rational curve given as $\pi_j^{-1} (\tilde{p}^{j-1})$
where $\tilde{p}^{j-1}$ is a singularity of $\wifol^{j-1}$.
Finally we identify curves and their proper transforms in the
obvious way. Also $\wiX^j$ will stand for the corresponding
blow-up of $Y$.

Throughout this section $Y, \fol$ are supposed to verify the
following assumptions:

\noindent {\bf A}. $Y = y^{-k} fZ$ where $k \geq 2$,
$Z$ is a holomorphic vector field having an isolated
singularity at $(0,0) \in \C^2$ and $f$ is a holomorphic
function.

\noindent {\bf B}. Assumptions~2 and~3 of Section~4.

\noindent {\bf C}. The origin is not a dicritical singularity
of $\fol$.

\noindent It immediately results from the above assumptions
that the axis $\{ y=0\}$ is a smooth separatrix of $\fol$.
Letting $Z= f\parlx + g \parly$, recall that the multiplicity of
$\fol$ along $\{ y=0\}$ (or the order of $\fol$ w.r.t.
$\{y=0\}$ and $(0,0)$) is by definition the order at $0
\in \C$ of the function $x \mapsto f(x,0)$.

The main result of this section is Theorem~(\ref{aque11}) below.

\begin{teo}
\label{aque11}
Let $Y, \fol$ be as above. Suppose that the divisor of zeros/poles
of $\wiX^r$ contains $\cale^r$ (i.e. there is no component $D^j$
of $\cale^r$ where $\wiX^r$ is regular). Then the multiplicity of
$\fol$ along $\{ y=0\}$ is at most~$2$ (in particular the
order of $\fol$ at $(0,0)$ is not greater than~$2$).
\end{teo}

As a by-product of our proof, the cases in which the multiplicity
of $\fol$ along $\{ y=0\}$ is~$2$ are going to be characterized as
well.
Also note that assumption~{\bf C} ensures that all the components
$D^j$ are invariant by $\fol^r$. Moreover none of the
singularities of $\wifol^r$ is dicritical. In what follows
we shall obtain the proof of Theorem~(\ref{aque11}) through
a systematic analyse of the structure of the resolution tree
of $\fol$.

\begin{obs}
\label{obs11}
{\rm Let $\fol$ be a foliation defined on a neighborhood of
$(0,0) \in \C^2$ and consider a separatrix $\cals$ of $\fol$.
Denote by $\wifol$ the blow-up of $\fol$ and by
$\widetilde{\cals}$ the proper transform of $\cals$. Naturally
$\widetilde{\cals}$ constitutes a separatrix for some singularity
$p$ of $\wifol$. In the sequel the elementary relation
\begin{equation}
{\rm Ind}_{p} (\wifol , \widetilde{\cals}) = {\rm Ind}_{(0,0)}
(\fol ,\cals) -1 \, \label{baixa1}
\end{equation}
will often be used.}
\end{obs}

Consider the resolution tree~(\ref{resotree}) of $\fol$.
By assumption $\cale^r$ contains a rational curve $D^r$
of self-intersection~$-1$. Blowing-down (collapsing) this
curve yields a foliation $\wifol^{r_1}$ together with the
total exceptional divisor $\cale^{r_1}$. If $\cale^{r_1}$
contains a rational curve with self-intersection~$-1$ where
all the singularities of $\wifol^{r_1}$ are simple, we then
continue by blowing-down this curve. Proceeding recurrently
in this way, we eventually arrive to a foliation $\wifol^{r_1}$,
$1 \leq r_1 < r$, together with an exceptional divisor
$\cale^{r_1}$ such that every irreducible component of $\cale^{r_1}$
with self-intersection~$-1$ contains a singularity of $\wifol^{r_1}$
with vanishing eigenvalues. Let $\wiX^{r_1}$ be the vector
field corresponding to $\wifol^{r_1}, \; \cale^{r_1}$,
using Proposition~(\ref{prop4.2})
we conclude the following:

\begin{lema}
\label{sequence1}
$\cale^{r_1}$ contains (at least) one rational curve $D^{r_1}$
of self-intersection~$-1$. Moreover if $p$ is a singularity
of $\wifol^{r_1}$ belonging to $D^{r_1}$ then either $p$ is simple
for $\wifol^{r_1}$ or $\wiX^{r_1}$ has one of the normal
forms $Z_{1,11},\;  Z_{0,12}, \; Z_{1,00}$ around $p$. Finally there
is at least one such singularity $p_1$ which is not simple
for $\wifol^{r_1}$.\qed
\end{lema}

The next step is to consider the following
description of the models $Z_{1,11},
\;  Z_{0,12}, \; Z_{1,00}$ (and their respective associated
foliations $\fol_{1,11} , \; \fol_{0,12}, \; \fol_{1,00}$)
which results at once from the
definition of these models given in Section~5. While
this description is slightly less precise than the previous
one, it emphasizes the properties more often used in the sequel.

\noindent $\bullet$ $Z_{1,11}, \; \fol_{1,11}$\,: $\fol_{1,11}$
has exactly 2 separatrizes $\cals_1 ,\cals_2$ which are
smooth, transverse and of index {\it zero}. $\fol_{1,11}$
has order~$2$ at the origin. The multiplicity of $\fol_{1,11}$
along $\cals_1 , \; \cals_2$ is~$2$.
The vector field $Z_{1,11}$
has poles of order~$1$ on each of the separatrizes $\cals_1$,
$\cals_2$.

\noindent $\bullet$ $Z_{0,12}, \; \fol_{0,12}$\,: $\fol_{0,12}$
has order~$2$ at the origin and 2 smooth, transverse separatrizes
coming from the separatrizes associated to the singularities of
eigenvalues~$-1,2$ of $\wifol_{0,12}$ (they are
denoted $\cals_1 ,\cals_2$ and called strong separatrizes of
$\fol_{0,12}$). The multiplicity of $\fol_{0,12}$
along $\cals_1 , \; \cals_2$ is~$2$ and the corresponding
indices are both~$-1$. $\fol_{0,12}$ has still a formal third
separatrix $\cals_3$, referred to as the weak separatrix of
$\fol_{0,12}$, which may or may not be convergent. The vector
field $Z_{0,12}$ has poles of order~$d \in \N^{\ast}$ on both $\cals_1 ,\cals_2$.

\noindent $\bullet$ $Z_{1,00}, \; \fol_{1,00}$\,: $\fol_{1,00}$
has order~$2$ and 1 smooth separatrix $\cals_1$ coming from the singularity
of $\wifol_{1,00}$ whose eigenvalues are~$-1,1$ which
will be called the strong separatrix of $\fol_{1,00}$.
Note that the multiplicity of $\fol_{1,00}$ along $\cals_1$
is~$2$ and the index of $\cals_1$ is {\it zero}. $\wifol_{1,00}$
also has 2 additional formal separatrizes $\cals_2 ,\cals_3$
coming from the weak invariant manifold of the saddle-nodes
singularities of $\wifol_{1,00}$ and, accordingly, called
the weak separatrizes of $\fol_{1,00}$. Naturally the weak
separatrizes of $\fol_{1,00}$ may or may not converge.
Finally the vector field $Z_{1,00}$ has poles of order~$d \in \N^{\ast}$
on $\cals_1$.

\begin{obs}
\label{combin}
{\rm The content of this remark will not be proved in these
notes and therefore will not be formally used either. Nonetheless
it greatly clarifies the structure of the combinatorial discussion
that follows. Consider a meromorphic vector field $X$ having a
smooth separatrix $\cals$. Using appropriate coordinates $(x,y)$
we can identify $\cals$ with $\{ y=0 \}$ and write $X$
as $y^d (f\parlx + yg \parly)$ where $d \in \Z$ and $f(x,0)$
is a non-trivial meromorphic function. We define the {\it asymptotic
order} of $X$ at $\cals$ (at $(0,0)$), ${\rm ord\,asy}_{(0,0)} (X ,\cals)$,
by means of the formula
\begin{equation} 
{\rm ord\,asy}_{(0,0)} (X, \cals) = {\rm ord}_0 f(x,0) + d.
{\rm Ind}_{(0,0)} (\fol ,\cals) \, \label{defioras}
\end{equation}
where $\fol$ is the foliation associated to $X$. It can be proved
that $0 \leq {\rm ord\,asy}_{(0,0)} (X ,\cals) \leq 2$ provided
that $X$ is semi-complete. Besides, if $\cals$ is induced by a global
rational curve (still denoted by $\cals$) and $p_1 ,\ldots,
p_s$ are the singularities of $\fol$ on $\cals$,
then we have
\begin{equation}
\sum_{i=1}^s {\rm ord\,asy}_{p_i} (X, \cals) = 2 \label{somou2}
\end{equation}
provided that $X$ is semi-complete on a neighborhood of $\cals$.
Note that the above formula indeed generalizes
formula~(\ref{2sr=2}).}
\end{obs}

Now we go back to the vector field $\wiX^{r_1}$ on a neighborhood
of $D^{r_1}$. Again we denote by $p_1 ,\ldots ,p_s$ the singularities
of $\wifol^{r_1}$ on $D^{r_1}$ and, without loss of generality,
assume that $\wiX^{r_1}$ admits one of the normal
forms $Z_{1,11},\;  Z_{0,12}, \; Z_{1,00}$ around $p_1$.

\begin{lema}
\label{sequence2}
All the singularities $p_2 ,\ldots ,p_s$ are simple for $\wifol^{r_1}$.
\end{lema}

\noindent {\it Proof}\,: We have to prove that it is not
possible to exist two singularities with one of the normal
forms $Z_{1,11},\;  Z_{0,12}, \; Z_{1,00}$ on $D^{r_1}$. Clearly
we cannot have two singularities of type $Z_{1,11}$
otherwise a ``generic'' regular orbit of $\wiX^{r_1}$ would contain
two singular points which
contradicts assumption~{\bf B}. Indeed the fact
that a ``generic'' regular orbit of $\wiX^{r_1}$ effectively
intersects the divisor of zeros of both singularities results
from the method employed in the preceding section. In the
present case the discussion is simplified since the singularity
appearing in the intersection of the two irreducible components
of the exceptional divisor is linear with eigenvalues~$-1,1$
(cf. description of $Z_{1,11}$).

To prove that the other
combinations are also impossible, it is enough to repeat the
argument employed in Section~5, in particular using the fact
that the order $d_1\neq 0$ of $\wiX^{r_1}$ on $D^{r_1}$ does
not depend of the singularity $p_i$. If the reader takes for
grant Formula~(\ref{somou2}), this verification can easily
be explained. In fact, the asymptotic order of $Z_{0,12}$
(resp. $Z_{1,00}$) with respect to its strong separatrizes
is already~$2$. Furthermore the asymptotic order of $Z_{1,11}$
with respect to its separatrizes is~$1$. Since the asymptotic
order of a semi-complete singularity cannot be negative, it
becomes obvious that two such singularities cannot co-exist
on $D^{r_1}$ provided that $Y^{r_1}$ is semi-complete.\qed

Now we analyse each of the three possible cases.

\noindent $\bullet$ The normal form of $\wiX^{r_1}$ around $p_1$
is $Z_{1,11}$: First recall that $D^{r_1}$ is an irreducible
component of order~$1$ of the pole divisor of $\wiX^{r_1}$.
Since the number of singularities of regular orbits of $\wiX^{r_1}$
is at most~$1$, it follows that none of the remaining singularities
$p_2 ,\ldots , p_s$ can be a LJ-singularity for $\wifol^{r_1}$.
Otherwise there would be another curve of {\it zeros} of
$\wiX^{r_1}$ which is not invariant by $\wifol^{r_1}$ so that ``generic''
regular orbits of $\wiX^{r_1}$ would have~$2$ singularities (cf. above).
By the same reason the holonomy of $D^{r_1} \setminus
\{ p_2 ,\ldots ,p_s\}$ with respect to $\wifol^{r_1}$ must
be trivial. This implies that none of the singularities
$p_2 ,\ldots ,p_s$ is a saddle-node for $\wifol^{r_1}$. In fact,
by virtue of Theorem~(\ref{selano}), a saddle-node must have
strong invariant manifold contained in $D^{r_1}$ which ensures
that the above mentioned holonomy is non-trivial. Then we conclude
that all the remaining singularities $p_2 ,\ldots ,p_s$ have
eigenvalues $m_i ,-n_i$ ($i=2,\ldots ,s$) with $m_i ,n_i \in
\N^{\ast}$. Once again the fact that the holonomy of
$D^{r_1} \setminus \{ p_2 ,\ldots ,p_s\}$ is trivial shows that
$m_i /n_i \in \N^{\ast}$. Finally Formula~(\ref{index2})
implies that $s=2$ and $m_2=n_2 =1$. An immediate application
of Formulas~(\ref{ordem2}) and~(\ref{equa2}) (or yet
Formula~(\ref{somou2})) shows that the separatrix of $p_2$ transverse
to $D^{r_1}$ is a component of order~$1$ of the pole divisor
of $\wiX^{r_1}$. Finally we denote by $Z_{1,11}^{(1)}$ (resp.
$\fol_{1,11}^{(1)}$) the local vector field (resp. holomorphic
foliation) resulting from the collapsing of $D^{r_1}$.
Summarizing one has:

\noindent $\bullet$ $Z_{1,11}^{(1)}, \; \fol_{1,11}^{(1)}$\,: The foliation
$\fol_{1,11}^{(1)}$ has exactly two separatrizes $\cals_1 ,\cals_2$
which are smooth, transverse and of indices respectively equal
to~$1$ and~$0$. The order of $\fol_{1,11}^{(1)}$ at the origin
is~$2$ as well as the multiplicity of $\fol_{1,11}^{(1)}$
along $\cals_1 , \cals_2$. The vector field $Z^{(1)}_{1,11}$
has poles of order~$1$ on $\cals_1 ,\cals_2$.

Now let us discuss the second case.

\noindent $\bullet$ The normal form of $\wiX^{r_1}$ around $p_1$
is $Z_{0,12}$\,: Recall that $D^{r_1}$ is an irreducible component
of order $d\neq 0$ of the pole divisor of $\wiX^{r_1}$. Repeating the
argument of the previous section, we see that the singularities
$p_2 , \ldots ,p_s$ can be neither saddle-nodes nor LJ-singularities.
Again this can directly be seen from Formula~(\ref{somou2}): by virtue
of Lemma~(\ref{le3.3}) and Theorem~(\ref{selano}), both
types of singularities in question have asymptotic order equal to~$1$.
Nonetheless the asymptotic order of $Z_{0,12}$ is already~$2$
which implies the claim. It follows that $\wifol^{r_1}$ has eigenvalues
$m_i, -n_i$ at each of the remaining singularities $p_2 ,
\ldots, p_s$ ($m_i ,n_i \in \N^{\ast}$). In
particular the index of $D^{r_1}$ w.r.t. $\wifol^{r_1}$ around each
$p_i$ is strictly negative. Hence Formula~(\ref{index2}) shows
that $p_1$ is, in fact, the unique singularity of $\wifol^{r_1}$ on
$D^{r_1}$. We denote by $Z_{0,12}^{(1)}$ (resp. $\fol_{0,12}^{(1)}$)
the local vector field (resp. holomorphic foliation) arising from
the collapsing of $D^{r_1}$. Thus:

\noindent $\bullet$ $Z_{0,12}^{(1)}, \; \fol_{0,12}^{(1)}$\,: The
order of $\fol_{0,12}^{(1)}$ at the origin is~$1$, besides the linear
part of $\fol_{0,12}^{(1)}$ is nilpotent. $\fol_{0,12}^{(1)}$ has
one strong separatrix $\cals_1$ obtained through the strong separatrix
of $\fol_{0,12}$ which is transverse to $D^{r_1}$. This separatrix
is smooth and has index {\it zero}, furthermore
the multiplicity of $\fol_{0,12}^{(1)}$ along $\cals_1$ is~$2$.
The foliation $\fol_{0,12}^{(1)}$ has still another formal weak
separatrix which may or may not converge.
Finally the vector field $Z_{0,12}^{(1)}$ has poles of order $d\neq 0$
on $\cals_1$.

Finally we have:

\noindent $\bullet$ The normal form of $\wiX^{r_1}$ around $p_1$
is $Z_{1,00}$\,:
Note that $D^{r_1}$ is an irreducible component
of order~$d\neq 0$ of the pole divisor of $\wiX^{r_1}$
(cf. description of the vector field $Z_{1,00}$).
As before
the remaining singularities cannot be saddle-nodes or LJ-singularities
(the asymptotic order of $Z_{1,00}$ w.r.t. $D^{r_1}$ is already~$2$).
It follows that $\wifol^{r_1}$ has eigenvalues $m_i, -n_i$
at the remaining singularities $p_2 ,\ldots ,p_s$ ($m_i, n_i
\in \N$).
Around each singularity $p_i$ ($i=2,\ldots ,s$), the vector
field $\wiX^{r_1}$ can be written as
$$
x_i^{-d} t_i^{k_i} h_i [m_i x_i (1+ {\rm h.o.t.}) \partial /\partial
x_i - n_i t_i (1+ {\rm h.o.t.}) \partial /\partial t_i]
$$
where $h_i (0,0) \neq 0$. We just have to repeat the argument of
Section~$5$, here we summarize the discussion by using the ``fact''
that the asymptotic order of $\wiX^{r_1}$ w.r.t. $D^{r_1}$ has
to be {\it zero} at each $p_i$. Indeed, this gives us that
$k_i =-1+dm_i /n_i$. Comparing this with Lemma~(\ref{le3.1}),
it results that $n_i=1$. Hence
Formula~(\ref{index2})
informs us that $s=2$ and $m_2 =n_2 =1$. It also follows that
$k_2 =d-1$. Let us denote by $Z_{1,00}^{(1)}$ (resp. $\fol_{1,00}^{(1)}$)
the local vector field (resp. holomorphic foliation) arising
from the collapsing of $D^{r_1}$. 

\noindent $\bullet$ $Z_{1,00}^{(1)}, \; \fol_{1,00}^{(1)}$\,: The
order of $\fol_{1,00}^{(1)}$ at the origin is~$2$ and it has
one strong separatrix $\cals_1$ obtained through the separatrix of $p_2$
which is transverse to $D^{r_1}$. This separatrix
is smooth and has index {\it zero}, furthermore
the multiplicity of $\fol_{1,00}^{(1)}$ along $\cals_1$ is~$2$.
The foliation $\fol_{1,00}^{(1)}$ has still two formal weak
separatrizes which may or may not converge.
Finally the vector field $Z_{1,00}^{(1)}$ has poles of order $d\neq 0$
on $\cals_1$.

Summarizing what precedes, we easily obtain the following analogue
of Lemma~(\ref{sequence1}):

\begin{lema}
\label{sequence4}
$\cale^{r_2}$ contains (at least) one rational curve $D^{r_2}$
of self-intersection~$-1$. Moreover if $p$ is a singularity
of $\wifol$ belonging to $D^{r_2}$ then either $p$ is simple
for $\wifol^{r_2}$ or $\wiX^{r_2}$ has one of the normal
forms $Z_{1,11},\;  Z_{0,12}, \; Z_{1,00}, \; Z_{1,11}^{(1)},
\; Z_{0,12}^{(1)}, \; Z_{1,00}^{(1)}$ around $p$. Finally there
is at least one such singularity $p_1$ which is not simple
for $\wifol^{r_2}$.\qed
\end{lema}

The argument is now by recurrence, we shall discuss only
the next step in details.
Again $p_1 ,\ldots ,p_s$ are the singularities
of $\wifol^{r_2}$ on $D^{r_2}$ and, without loss of generality,
$\wiX^{r_2}$ admits one of the normal
forms $Z_{1,11},\;  Z_{0,12}, \; Z_{1,00}, \;Z_{1,11}^{(1)},
\; Z_{0,12}^{(1)}, \; Z_{1,00}^{(1)}$ around $p_1$.

\begin{lema}
\label{sequence5}
All the singularities $p_2 ,\ldots ,p_s$ are simple for $\wifol^{r_2}$.
\end{lema}

\noindent {\it Proof}\,: The proof is as in Lemma~(\ref{sequence2}).
Since no leaf of $\wifol^{r_2}$ can meet the divisor of zeros
of $\wiX^{r_2}$ in more than one point, it results that we can
have at most one singularity $Z_{1,11}$ or $Z_{1,11}^{(1)}$ on
$D^{r_2}$. The fact that the models $Z_{0,12}, \; Z_{1,00}, \;
Z_{0,12}^{(1)}$ and $Z_{1,00}^{(1)}$ cannot be combined among them or with
$Z_{1,11}, \; Z_{1,11}^{(1)}$ follows from the natural
generalization of the method of Section~5 (which is again
explained by the fact that these $3$ vector fields have
asymptotic order~$2$ w.r.t. $D^{r_2}$). The lemma is proved.\qed

So we have obtained three new possibilities according to the normal
form of $\wiX^{r_2}$ around $p_1$ is $Z_{1,11}^{(1)}, \;
Z_{0,12}^{(1)}$ or $Z_{1,00}^{(1)}$. Let us analyse them separately.

\noindent $\bullet$ The normal form of $\wiX^{r_2}$ around $p_1$
is $Z_{1,11}^{(1)}$\,: Note that $\fol_{1,11}^{(1)}$ has two
separatrizes which may coincide with $D^{r_2}$, one with index
{\it zero} and other with index~$1$. In any case, $D^{r_2}$
is an irreducible component of order~$1$ of the divisor of poles
of $\wiX^{r_2}$. Suppose first that the index of $D^{r_2}$ w.r.t.
$\fol^{r_2}$ at $p_1$ is {\it zero}. The discussion then goes
exactly as in the case of $Z_{1,11}$. We conclude that $s=2$
and that $\fol^{r_2}$ has eigenvalues~$-1,1$ at $p_2$. The fact
that the holonomy of $D^{r_2} \setminus \{ p_1 ,p_2 \}$ w.r.t.
$\fol^{r_2}$ is trivial also implies that $\fol^{r_2}$ is
linearizable at $p_2$. Formulas~(\ref{ordem2}) and~(\ref{equa2})
show that the separatrix of $p_2$ transverse to $D^{r_2}$ is
a component with order~$1$ of the divisor of poles of $\wiX^{r_2}$.
Finally we denote by $Z_{1,11}^{(2)}$ (resp. $\fol_{1,11}^{(2)}$)
the local vector field (resp. holomorphic foliation) arising from
the collapsing of $D^{r_2}$. One has

\noindent $\bullet$ $Z_{1,11}^{(2)}, \; \fol_{1,11}^{(2)}$\,:
The foliation $\fol_{1,11}^{(2)}$
has exactly 2 separatrizes $\cals_1 ,\cals_2$ which are
smooth, transverse and of indices respectively equal
to {\it zero} and~$2$. $\fol_{1,11}^{(2)}$
has order~$2$ at the origin and its multiplicity
along $\cals_1 , \; \cals_2$ is~$2$.
The vector field $Z_{1,11}^{(2)}$
has poles of order~$1$ on each of the separatrizes $\cals_1$,
$\cals_2$.

Now let us prove that the index of $D^{r_2}$ w.r.t. $\fol^{r_2}$
at $p_1$ cannot be~$1$. Suppose for a contradiction that this
index is~$1$. Again the triviality of the holonomy of the regular
leaf contained in $D^{r_2}$ implies that all the singularities
$p_2 ,\ldots ,p_s$ are linearizable with eigenvalues~$-1,1$.
Hence Formula~(\ref{index2}) ensures that $s=3$. In turn,
Formula~(\ref{equa2}) shows that the separatrix of $p_2$
(resp. $p_3$) transverse to $D^{r_2}$ is a component with
order~$1$ of the pole divisor of $\wiX^{r_2}$. This is however
impossible in view of Formula~(\ref{somou2}). Alternate, we
can observe that the divisor of poles
of the vector field obtained by collapsing $D^{r_2}$ consists
of three smooth separatrizes. By virtue of assumptions {\bf A,
B, C} one of them must be the proper transform of $\{ y=0\}$.
In particular its order as component of the pole divisor
should be $k \geq 2$, thus providing the desired contradiction.

Next one has:

\noindent $\bullet$ The normal form of $\wiX^{r_2}$ around $p_1$
is $Z_{0,12}^{(1)}$\,: The divisor $D^{r_2}$ constitutes a separatrix
of $\wifol^{r_2}$ at $p_1$.
Besides $\wiX^{r_2}$ has poles of order $d \neq 0$ on $D^{r_2}$
(cf. description
of $Z_{0,12}^{(1)}$). Note also that the index of
$D^{r_2}$ w.r.t. $\wifol^{r_2}$ at $p_1$ is {\it zero}.
Our standard method shows that
the remaining singularities cannot be saddle-nodes or LJ-singularities
(the asymptotic order of $Z_{0,12}^{(1)}$
w.r.t. $D^{r_2}$ is already~$2$). Thus $\wifol^{r_2}$ has eigenvalues $m_i, -n_i$
at  the singularity $p_i$, $i=2,\ldots s$ ($m_i, n_i
\in \N$). On a neighborhood of $p_i$, $\wiX^{r_2}$ is
given in appropriate coordinates by
$$
x_i^{-d} t_i^{k_i} h_i [m_i x_i (1+ {\rm h.o.t.}) \partial /\partial
x_i - n_i t_i (1+ {\rm h.o.t.}) \partial /\partial t_i]
$$
where $h_i (0,0) \neq 0$. It is enough to repeat the discussion of
Section~$5$. Using the ``fact''
that the asymptotic order of $\wiX^{r_1}$ w.r.t. $D^{r_1}$ has
to be {\it zero} at each $p_i$, this can be summarized as follows.
The asymptotic order is given by $k_i + 1 -dm_i /n_i$, since it
must equal {\it zero}, one has $k_i =-1 +dm_i/n_i$. Comparing
with Lemma~(\ref{le3.1}), we conclude that $n_i=1$. Hence
Formula~(\ref{index2}) implies that $s=2$ and $m_2=n_2 =1$.
In particular $k_2=d-1$. We then denote by $Z_{0,12}^{(2)}$
(resp. $\fol_{0,12}^{(2)}$) the local vector field
(resp. holomorphic foliation) arising from the collapsing of
$D^{r_2}$.

\noindent $\bullet$ $Z_{0,12}^{(2)}, \; \fol_{0,12}^{(2)}$\,: The
order of $\fol_{0,12}^{(2)}$ at the origin is~$2$. $\fol_{0,12}^{(2)}$ has
one strong separatrix $\cals_1$ obtained through the strong separatrix
of $\fol_{0,12}^{(1)}$ which is transverse to $D^{r_2}$. This separatrix
is smooth and has index {\it zero}, furthermore
the multiplicity of $\fol_{0,12}^{(2)}$ along $\cals_1$ is~$2$.
The foliation $\fol_{0,12}^{(2)}$ has still another formal weak
separatrix which may or may not converge.
Finally the vector field $Z_{0,12}^{(2)}$ has poles of order $d\neq 0$
on $\cals_1$.

Finally let us consider $Z_{1,00}^{(1)}$.

\noindent $\bullet$ The normal form of $\wiX^{r_2}$ around $p_1$
is $Z_{1,00}^{(1)}$\,: The discussion is totally analoguous to
the case $Z_{1,00}$. After collapsing $D^{r_2}$, we obtain a
local vector field $Z_{1,00}^{(2)}$ (resp. holomorphic
foliation $\fol_{1,00}^{(2)}$) with the following characteristics:

\noindent $\bullet$ $Z_{1,00}^{(2)}, \; \fol_{1,00}^{(2)}$\,: The
order of $\fol_{1,00}^{(2)}$ at the origin is~$2$ and it has
one strong separatrix $\cals_1$ obtained through the separatrix of $p_2$
which is transverse to $D^{r_2}$. This separatrix
is smooth and has index {\it zero}, furthermore
the multiplicity of $\fol_{1,00}^{(2)}$ along $\cals_1$ is~$2$.
The foliation $\fol_{1,00}^{(2)}$ has still two formal weak
separatrizes which may or may not converge.
Finally the vector field $Z_{1,00}^{(2)}$ has poles of order $d\neq 0$
on $\cals_1$.

Let us inductively define a sequence of vector fields $Z_{1,11}
^{(n)}$ by combining over a rational curve with self-intersection~$-1$
a model $Z_{1,11}^{(n-1)}$ with a linear singularity $p_2$
having eigenvalues~$-1,1$. The rational curve in question induces
a separatrix of index {\it zero} for $Z_{1,11}^{(n-1)}$ and both
separatrizes of $p_2$ are components having order~$1$ of the
pole divisor of the corresponding vector field. The model
$Z_{1,11}^{(n)}$ is then obtained by collapsing the
mentioned rational curve. Similarly we also define the sequences
$Z_{1,00}^{(n)}, \, Z_{0,12}^{(n)}$.

\vspace{0.1cm}

\noindent {\it Proof of Theorem~(\ref{aque11})}\,: Let $Y,\;
\fol$ be as in the statement of this theorem. Suppose
first that the order of $\fol$ at $(0,0)$ is greater than one.
We consider a resolution tree~(\ref{resotree}) for $\fol$ and
the recurrent procedure discussed above. Whenever we collapse
a rational curve with self-intersection~$-1$ contained in
one of the exceptional divisors $\cale^j$, it results a singularity
which is either simple or belongs to the list $Z_{1,11},
\; Z_{1,11}^{(n)}, \; Z_{1,00}, \; Z_{1,00}^{(n)},\; Z_{0,12}, \;
Z_{0,12}^{(n)}$.
After a finite number of steps of the above procedure, we arrive
to the original vector field $Y$ (resp. foliation $\fol$). Therefore
$Y$ must admit one of the above indicated normal forms. Since the
divisor of poles of $Y$ is constituted by the axis $\{ y=0\}$,
the cases $Z_{1,11},\; Z_{1,11}^{(n)},  \; Z_{0,12}$ cannot be
produced (their divisor of poles consist of two irreducible components).
The case $Z_{0,12}^{(1)}$ cannot be produced either since the order
of the associated foliation is supposed to be greater than or
equal to~$2$. Thus we conclude that $Y=Z_{1,00}^{(n)}$ or,
for $n\geq 2$, $Y=Z_{0,12}^{(n)}$ and the theorem
follows in this case.

Next suppose that the order of $\fol$ at $(0,0)$ is~$1$.
Clearly if the linear part of $\fol$ at $(0,0)$ has rank~$2$
(i.e. the corresponding holomorphic vector field with isolated
singularities has~$2$ non-vanishing eigenvalues at $(0,0)$),
then the conclusion is obvious. Next suppose
that $\fol$ is a saddle-node. If $\{ y=0\}$ corresponds to
the strong invariant manifold of $\fol$ then the statement
is obvious. On the other hand, $\{ y=0\}$ cannot be the
weak invariant manifold thanks to Theorem~(\ref{selano}).

It only remains to check the case where the linear part
of $\fol$ at $(0,0)$ is nilpotent. The blow-up
$\wiX$ (resp. $\wifol$) of $Y$ (resp. $\fol$) is such that
$\wifol$ has a single singularity $p \in \pi^{-1} (0)$. In
addition the order of $\wifol$ at $p$ is necessarily~$2$. On
a neighborhood of $p$, the vector field $\wiX$ is given by
$x^{-k} y^{-k} Z$ where $Z$ is as in assumption~{\bf A}. The
assumptions~{\bf B}, {\bf C} are clearly verified as well.
An inspection in the preceding discussion immediately shows
that it applies equally well to this vector field $\wiX$.
We conclude that $\wiX$ is given on a neighborhood of $p$
by the model $Z_{0,12}$. Hence $Y$ is the model
$Z_{0,12}^{(1)}$ completing the proof of the theorem.\qed

\bigskip

\noindent {\bf {\large $\bullet$ Conclusion}}:

\vspace{0.1cm}

\noindent {\it Proof of Theorem~A}\,: Let $X$ be a complete polynomial
vector field in $\C^2$ with degree~$2$ or greater.
We denote by $\fol_X$ the singular holomorphic
foliation induced by $X$ on $\pp$. We know from Lemma~(\ref{lema2.1})
that the line at infinity $\Delta \subset \pp$ is invariant
under $\fol_X$. On the other hand there is a dicritical
singularity $p_1$ of $\fol_X$ belonging to $\Delta$. Hence
Proposition~(\ref{jouano}) ensures that $\fol_X$ has a meromorphic
first integral on $\pp$. Besides the generic leaves of $\fol_X$
in $\pp$ are, up to normalization, rational curves (i.e. isomorphic
to $\C \Pp (1)$). According to Saito and Suzuki (cf. \cite{suzu2}),
up to polynomial automorphisms of $\C^2$, $\fol_X$ is given by
a first integral $R$ having one of the following forms:
\begin{description}

\item[$\imath$)] $R(x,y) = x$;

\item[$\imath \imath$)] $R(x,y) = x^n y^m$, with
${\rm g.c.d} (m,n) =1$ and $m,n \in \Z$;

\item[$\imath \imath \imath$)] $R(x,y)=x^n(x^ly + P(x))^m$,
with ${\rm g.c.d} (m,n) =1$ and $m,n \in \Z$, $l\geq 1$.
Moreover $P$ is a polynomial of degree at most $l-1$ satisfying
$P (0) \neq 0$.
\end{description}
To each of these first integrals there corresponds the foliations
associated to the vector fields $X_1 =\parlx$,
$X_2 = mx \parlx - ny \parly$ and $X_3 = mx^{l+1} \parlx
-[(n+lm)x^ly + nP(x) + mxP' (x)] \parly$.
Therefore the original vector field $X$ has the form
$X=Q. X_i$, where $Q$ is a polynomial and $i =1,2,3$.
If $i=1$ then it follows at once that $Q$ has to be as
in the in order to produce a complete vector field $X$.
Assume now that $i=2$. Using for instance Lemma~(\ref{le3.1}),
we see that $P$ has again the form indicated in the statement
unless $X_2 = x \parlx - y \parly$ in which case we can also
have $P = (xy)^a (x-y)$. Nonetheless we immediately check that
the resulting vector field $X$ is not complete in this case.

Finally let us assume that $i=3$. It is again easy to see that
that the resulting vector field cannot be complete. This follows
for example from the fact that $(0,0)$ is a singularity of $X_3$
with trivial eigenvalues (cf. \cite{re4}). The theorem is proved.\qed

\bigskip

\noindent {\it Proof of Theorem~B}\,: Let us suppose for a contradiction
that none of the singularities $p-1 ,\ldots ,p_k$ of $\fol_X$
in $\Delta$ is dicritical. We write $X$ as $F.Z$ where $F$ is
a polynomial of degree $n\in \N$ and $Z$ is a polynomial vector
field of degree $d-n$ and having isolated zeros (where $d$
is the degree of $X$).

We consider the restriction of $X$ to a neighborhood of $p_i$.
Clearly this restriction satisfies assumptions~{\bf A}, {\bf B}
and~{\bf C} of Section~5. In particular Theorem~(\ref{aque11})
applies to show that the multiplicity of $\fol$ along $\Delta$
is at most~$2$. Moreover, if this multiplicity is~$2$, then
$X$ admits one of the normal forms
$Z_{1,11},
\; Z_{1,11}^{(n)}, \; Z_{1,00}, \; Z_{1,00}^{(n)},\; Z_{0,12}, \;
Z_{0,12}^{(n)}$ on a neighborhood of $p_i$.

From Lemma~(\ref{lema2.4}) we know that $k \leq 3$.
Since the sum of the multiplicities of $\fol$ along $\Delta$
at each $p_i$ is equal to $d-n+1$, it follows that $d-n \leq 5$.

However, if $k=3$, Corollary~(\ref{lema2.3}) shows that the
top-degree homogeneous component of $X$ is as in the cases~5, 6, 7
of the corollary in question. Simple calculations guarantees
that it is not possible to realize the models $Z_{1,11},
\ldots , Z_{0,12}^{(n)}$ in this way. The case $k=1$ being
trivial, we just need to consider the case $k=2$. Now we
have $d-n \leq 3$ and again it is very easy to conclude
that none of these possibilities lead to a complete
polynomial vector field. The resulting contradiction proves
the theorem.\qed

\noindent {\sc Julio C. Rebelo
\begin{tabbing}
{\rm {\bf Permanent Address}}  \hspace{3.8cm}   \=  {\rm {\bf Current Address}}\\
PUC-Rio                    \>    Institute for Mathematical Sciences\\
R. Marques de S. Vicente 225 \> SUNY at Stony Brook\\
Rio de Janeiro RJ CEP 22453-900 \>         Stony Brook NY 11794-3660\\
Brazil                   \>             USA\\
{\it email jrebelo@mat.puc-rio.br} \>      {\it email jrebelo@math.sunysb.edu}
\end{tabbing}

}

\end{document}